\documentclass[12pt]{article}

\usepackage{amsgen,amsmath,amstext,amsbsy,amsopn,amsfonts,amssymb}
\usepackage[dvips]{graphicx}

\title{A non radical based approach to study of associative algebras}
\author{Vladimir Dergachev}

\def\buzz#1{{\normalfont{\bf{[#1]}}}}

\def\ind{\textrm{\normalfont ind\,}}

\def\ad#1{\textrm{\normalfont ad}_{{#1}}\,}
\def\Ad#1{\textrm{\normalfont Ad}_{{#1}}\,}

\def\coAd#1{\textrm{\normalfont Ad}^*_{{#1}}\,}
\def\diag{\textrm{\normalfont diag}}
\def\Mat{\textrm{\normalfont Mat}}

\def\Stab{\textrm{\normalfont Stab}}

\def\Nil{\textrm{\normalfont Nil}}

\def\ker{\textrm{\normalfont ker}\,}
\def\AF{A_F}
\def\AFVA#1{A_{{#1}}}
\def\AFVAVA#1#2{{{#1}}_{{#2}}}

\def\imply{\Rightarrow}
\def\implies{\imply}
\def\follows{\imply}
\def\intersect{\cap}

\def\image{\textrm{\normalfont Im}\,}
\def\span{\textrm{\normalfont span\,}}
\def\codim{\textrm{\normalfont codim\,}}
\def\Spec{\textrm{\normalfont Spec\,}}

\newtheorem{theorem}{Theorem}
\newtheorem{definition}{Definition}
\newtheorem{example}{Example}
\newtheorem{proposition}[theorem]{Proposition}

\newtheorem{lemma}[theorem]{Lemma}

\newenvironment{proof}{\begin{trivlist}\item[\hskip%
\labelsep{\bf Proof\quad}]}%
{\hfill\qed\end{trivlist}}
\newcommand{\qed}{{\unskip\nobreak\hfil\penalty50\hskip .001pt \hbox{}
          \nobreak\hfil
          \vrule height 1.2ex width 1.1ex depth -.1ex
           \parfillskip=0pt\finalhyphendemerits=0\medbreak}\rm}

\begin{document}
\maketitle
\abstract{
We  study pairs of associative algebras and linear functionals.
New results together with corrected proofs for previously published
material are presented. In particular, we prove the identity 
$\ind\Mat_n \otimes \mathfrak A=n \cdot\ind \mathfrak A$ for finite-dimensional
unital associative algebra $\mathfrak A$ with index~$1$.
}
\tableofcontents

\section{Introduction}
In this paper we describe our progress in the study of interaction between 
associative algebras and linear functionals defined on them. This aspect of associative
algebras is related to the following classical concepts:
\begin{itemize}
\item The Orbit method in the theory of Lie algebras. As any associative algebra
can be converted into a Lie algebra one hopes the additional
structure present in associative algebras will expose new phenomena.
\item The classical notion of multiplicative functionals. As will be shown later
multiplicative functionals are exactly the functionals whose
associated bilinear form has rank $1$.
\item Hopf algebras. The pair $(\textrm{associative algebra}, \textrm{functional})$ can be considered as intermediate concept between
associative algebras and Hopf algebras.
\end{itemize}

This investigation has been prompted by the observation (see \cite{DK}) that for a class of subalgebras of matrix algebra
the index in Lie algebra sense (i.e. the dimension of the kernel of Kirillov's form $B_F$ in generic functional $F$) 
of a tensor product of $\Mat_n$ with the algebra $\mathfrak A$ was exactly $n$ (which is the index of $\Mat_n$) times the index of the algebra $\mathfrak A$.

This identity does not readily generalize to an arbitrary Lie algebra nor to an arbitrary pair of associative algebras. The last theorem
in this paper establishes that the identity does hold for two type $1$ associative algebras provided that the pair satisfies some additional conditions.
These conditions are fulfilled by $\Mat_n$ and unital associative algebra of index $1$.

The question of expressing index of a tensor product of two associative algebras via known invariants of these algebras is still open, 
even if one restricts consideration to type $1$ algebras.

The method developed for proving this result has several interesting properties:
\begin{itemize}
\item the decomposition obtained can be considered an exponentiated version of root spaces decomposition of Lie algebras (and is, in fact,
exactly so for $\Mat_n$, with Cartan subalgebra being $\Stab_F(1)$). This decomposition is defined for any
associative algebra and has proved very convenient in analyzing coadjoint representation.
\item besides the Kirillov's form defined on coadjoint representation of associative algebras
one obtains a quadratic form on the stabilizer of coadjoint action. The non-degeneracy of the latter
corresponds to type $1$ algebras. For regular functionals, non-degeneracy of this quadratic form
implies that $\Stab_F(1)$ is a Frobenius algebra.
\item characteristic polynomial presents an easy way to obtain invariants of coadjoint action.
\end{itemize}

Some of the material has appeared in our earlier preprints \cite{D}, \cite{DD} and \cite{DDD} that mark
the progress of our study. Besides corrections to proofs, this paper
refines notation for characteristic polynomial and spaces $\Stab_F(\alpha)$ and $V(\alpha)$. We also introduce
the definitions of $\alpha(F)$-regular and $\alpha$-precise functionals.

Also new are the definitions of three types of algebras. The most studied ones are type $1$ with many of results
having generalizations to type $2$ algebras. 

Type $3$ algebras do not allow decomposition into direct sum
of spaces $V(\alpha)$ even after factoring by $\Nil_F$, but possess an interesting property of having 
non-empty $\Stab_F(\alpha)$ with non-trivial multiplication table for any $\alpha$. This awaits further study.

\section{Definitions, characteristic polynomial}
Before proceeding further let us introduce some definitions. 

Let $\mathfrak A$ be an associative algebra. Unless specially noted we will assume
$\mathfrak A$ to be a finite-dimensional algebra over complex field.
We denote by $\mathfrak A^*$ the dual space of linear functionals on $\mathfrak A$.

The multiplication law $\mathfrak A\otimes\mathfrak A\rightarrow \mathfrak A$ can
be considered as an $\mathfrak A$-valued bilinear form $A$ on $\mathfrak A$.
If one picks a basis $\left\{e_k\right\}$ in the algebra $\mathfrak A$ then $A$
can be represented by a $\mathfrak A$-valued matrix $\left(e_ie_j\right)$ in
this basis. Usually we will abuse notation by denoting this matrix by the same
letter $A$.

For a linear functional $F\in\mathfrak A^*$ we denote by $\AF$ the
bilinear form $F(A(\cdot,\cdot))$.

We denote by $\Stab_F(0)=\ker \AF$ and by $\Stab_F(\infty)=\ker \AF^T$.
Let $\Nil_F=\Stab_F(0)\intersect\Stab_F(\infty)$.

We also introduce the characteristic polynomial 
$$
\chi_{V}(\lambda,\mu, F)=\det_V\left(\lambda\AF+\mu\AF^T\right)
$$
where determinant is evaluated in some basis of vector space $V\subset \mathfrak A$ - this polynomial
is thus defined up to a constant multiple. 

In places where the functional $F$ is fixed we will use the notation $\chi_{F,V}(\lambda, \mu)$ - this saves a little
space in formulas. Lastly, we omit $V$ when $V=\mathfrak A$.

We will distinguish three configurations $(\mathfrak A, F)$:
\begin{itemize}
\item[\bf Type 1:] The characteristic polynomial of the entire algebra $\chi_{F,\mathfrak A}$
does not vanish. This implies $\dim \Nil_F=0$.
\item[\bf Type 2:] $\Nil_F$ has positive dimension and the characteristic polynomial of the subspace of maximal dimension
that is transversal to $\Nil_F$ does not vanish.
\item[\bf Type 3:] The characteristic polynomial of the vector space of maximal
dimension that is transversal to $\Nil_F$ does vanish.
\end{itemize}

We will consider an associative algebra $\mathfrak A$ to be of type $N$ if its dual space
$\mathfrak A^*$ possesses non-empty open subset (in either Zariski or Euclid topology) 
of functionals $F$ that form a pair of type $N$ with $\mathfrak A$.

We will now establish correctness of these definitions.

\begin{definition} Fix a certain topology on the space of linear functionals on $\mathfrak A$.
We will call a condition on linear functionals $F$ generic if, for a given associative algebra,
it is either not satisfied for any functional or is satisfied for an open dense set of functionals.
\end{definition}

\begin{theorem} Fix a subspace $V$ inside each finite-dimensional associative algebra over $\mathbb C$. Then the condition
that $\chi_{F,V}(\lambda, \mu)$ does not vanish as polynomial in $\lambda$ and $\mu$ is generic in both
Euclid and Zariski topologies.
\end{theorem}
\begin{proof}
First of all let us note that for a fixed subspace $V\in\mathfrak A$ the existence of non-empty set of functionals $F$
for which $\chi_{F,V}(\lambda, \mu)$ does not vanish is equivalent to non-vanishing
of the polynomial $\chi_V(\lambda, \mu, F)$ in all three variables. 

Secondly, the condition that $\chi_{F,V}(\lambda, \mu)$ vanishes
can be written as a system of polynomial equations in $F$ by equating coefficients at $\lambda$ and $\mu$ in $\chi_V(\lambda, \mu, F)$
to zero. If there is a point $F$ when at least one of these coefficients is non-zero then, by continuity, there exists an open neighbourhood
(in either Zariski or Euclid topology) in which this coefficient does not vanish and thus $\chi_{F,V}(\lambda, \mu)$
does not vanish at all in this neighbourhood. 

Lastly if there were an open neigbourhood
such that for all $F$ in it the coefficients vanish this would imply that the coefficients vanish identically in $F$ as they are polynomial.

Thus either the set of functionals $F$ for which $\chi_{F,V}(\lambda, \mu)$ does not vanish is empty or it is open and dense.
\end{proof}

\begin{theorem} The condition that a functional $F$ has the smallest $\dim \Nil_F$ is generic
in both Euclid and Zariski topologies.
\end{theorem}
\begin{proof}
Indeed, consider any functional $F_0$ that has the property that $\dim \Nil_{F_0}$ 
is minimal (it exists as the $\dim\Nil_F$ is a non-negative integer).
The coefficients of the linear system that defines $\Nil_F$ 
are themselves linear in $F$. Therefore the minors of this system are polynomials in $F$. There must be a minor of 
dimension $\dim \mathfrak A -\dim \Nil_{F_0}$ that does not vanish for a dense open set (in either Zariski or Euclid topologies)
of functionals $F$ that includes $F_0$. However, all of these functionals must have $\dim \Nil_F$ not less than $\dim \Nil_{F_0}$. 
Therefore, the set of all functionals $F$ with minimal $\dim \Nil_F$ is a union of dense open sets and thus is itself
open and dense.

This condition is always satisfied by some functionals as $\dim \Nil_{F_0}$ is a non-negative integer.
\end{proof}

\begin{theorem}A finite-dimensional associative algebra over $\mathbb C$ is either type $1$, $2$
or $3$. The definition does not change whether one considers Zariski or Euclid topologies.
\end{theorem}
\begin{proof}

As proved above, the condition that $\chi_{F,\mathfrak A}(\lambda, \mu)$ does not vanish identically is generic and
implies $\dim \Nil_F=0$. Therefore, type $1$ algebras are mutually exclusive with type $2$ or $3$.

Consider now the the case of associative algebra with $\dim \Nil_F=0$. Either $\chi_{\mathfrak A}(\lambda, \mu, V)$
does not vanish, in which case it is type $1$, or it does vanish, in which case it is type $3$. Thus, 
finite-dimensional associative algebra with $\dim \Nil_F=0$ is either type $1$ or type $3$.

Now we will concentrate on the situation where the minimal dimension of $\Nil_F$ is positive.

For each associative algebra pick $F_0$ such that $\dim\Nil_{F_0}$ is minimal and pick a subspace $V_0$ of dimension
$\dim \mathfrak A -\dim \Nil_{F_0}$ that is transversal to $\Nil_{F_0}$.

We know that the condition that $\chi_{F, V_0}(\lambda, \mu)$ does not vanish is generic. Therefore, for each
associative algebra,  either there
exists a dense open set of functionals which possess a subspace of maximal dimension (i.e. $V_0$) with non-vanishing characteristic
polynomial or there is a dense open set of functionals which possess a subspace of maximal dimension (i.e. $V_0$) on which
characteristic polynomial vanishes.

We will now prove that, for a given functional $F$, the characteristic polynomial $\chi_{F, V}(\lambda, \mu)$ 
either vanishes for all subspaces $V$ of maximal dimension that are transversal to $\dim \Nil_F$ or does not
vanish for any such $V$.

Consider a basis of $\mathfrak A$ subordinate to the direct sum $\mathfrak A=\Nil_{F}\oplus V$.
In this basis linear automorphisms of $\mathfrak A$ (as a vector space) that preserve
$\Nil_{F}$ have the following block structure:
$$
\begin{array}{c|cc}
 & \Nil_{F} & V \\
 \hline 
\Nil_{F} &  T_{NN} & 0 \\
V  &  T_{NV} & T_{VV} \\
\end{array}
$$

These automorphisms act transitively on the set of maximal subspaces $V$ that are transversal to $\Nil_F$.

The multiplication table written in this basis  has zeros in 
all rows and columns corresponding to basis vectors from $\Nil_{F}$. When acted on
by linear transformation that preserved $\Nil_{F}$ the multiplication table will still
have zeros in all rows and columns corresponding to basis vectors from $\Nil_{F}$. Furthermore,
the entries corresponding to two basis vectors from $V$ will only depend on $T_{VV}$ - an inner
automorphism of $V$. Therefore the property that the minor formed by restriction of $\AF$ to $V$
is zero or not is independent of the choice of subspace $V$.

Thus, an algebra with positive minimal $\dim \Nil_F$ either possesses a dense open set of functionals $F$
that have a maximal subspace transversal to $\Nil_F$ with non-vanishing characteristic polynomial (and so it is
type $2$) or there is a dense open set of functionals $F$ (which includes those with minimal $\dim \Nil_F$)
for which the characteristic polynomial vanishes on any maximal subspace transversal to $\Nil_F$ (and so it is type $3$).
\end{proof}

\section{Examples}
We will now present examples of associative algebras of all three types.

{\bf We will use the following notation: when writing multiplication tables letters denote basis elements in
$\mathfrak A$ and when writing characteristic polynomial we will use the same letters to denote value
of generic functional $F$ on this element. 

In other words, we are considering $\mathfrak A$ as linear functions
on $\mathfrak A^*$ and we compute characteristic polynomials by using multiplication of $S(\mathfrak A)$ (i.e. we multiply them
as polynomials over $\mathfrak A^*$), not the
multiplication of the associative algebra itself. }
\subsection{Type 1}
\begin{example}\normalfont \buzz{$\Mat_2$}
Let $a$,$b$,$c$,$d$ denote the 
matrix units $E_{1,1}$,$E_{1,2}$,$E_{2,1}$ and
$E_{2,2}$ correspondingly. Then the multiplication table $A$ is 
$$
\begin{array}{cccccc}
&\vline & a & b & c & d \\
 \hline
a& \vline& a & b & 0  & 0 \\
b& \vline& 0 & 0 & a & b  \\
c&\vline& c  & d  & 0 & 0 \\
d& \vline& 0 & 0  & c   & d \\
\end{array}
$$

The characteristic polynomial is equal to
$$
\begin{array}{l}
\chi(\lambda,\mu,F)=\det(\lambda A+\mu A^T)=\\
\qquad\qquad=-(\lambda+\mu)^2(ad-bc)(
(\lambda-\mu)^2(ad-bc)+\lambda\mu(a+d)^2))
\end{array}
$$
There are plenty of functionals $F$ for which the above expression does not vanish.
\end{example}

\begin{example}\normalfont\buzz{Seaweed 12x21} Let $\mathfrak A$ be the following 
subalgebra of $\Mat_3(\mathbb C)$:
$$
\left(
\begin{matrix}
a & b & 0 \\
0 & c & 0 \\
0 & d & e
\end{matrix}
\right)
$$
The multiplication table $A$ of $\mathfrak A$ is 
$$
\begin{array}{ccccccc}
&\vline & a & b & c & d & e \\
\hline 
a & \vline & a & b & 0 & 0 & 0\\
b & \vline & 0 & 0 & b & 0 & 0\\
c & \vline & 0 & 0 & c & 0 & 0\\
d & \vline & 0 & 0 & d & 0 & 0\\
e & \vline & 0 & 0 & 0 & d & e
\end{array}
$$
The characteristic polynomial is equal to
$$
\chi(\lambda,\mu,F)=\lambda^2\mu^2(\lambda+\mu)b^2d^2(a+c+e)
$$
As in the previous example the set of functionals $F$ for which characteristic polynomial
of the entire algebra does not vanish is Zariski open.
\end{example}
\label{mat_n}
\begin{example}\normalfont\buzz{$\Mat_n$}\label{Mat_n}
\begin{theorem}\label{chi_invariance}
The characteristic polynomial of the entire algebra is quasi-invariant under coadjoint action.
That is 
$$
\chi\left(\lambda,\mu,\coAd{g}F\right)=\left(\det \Ad{g}\right)^{-2}\chi\left(\lambda,\mu,F\right)
$$
\end{theorem}
\begin{proof}
Indeed, the matrix element $(i,j)$ of $\AF$ is
given by the expression $F(e_ie_j)$. Since
$$
\left(\coAd{g}F\right)(e_ie_j)=F(g^{-1}e_ie_jg)=F((g^{-1}e_ig)(g^{-1}e_jg))
$$
the substitution $F\rightarrow \coAd{g}F$ is  equivalent to the change of
basis induced by the matrix $\textrm{\normalfont{Ad}}_{g}^{-1}$. 
\end{proof}
\begin{definition}\buzz{Generalized resultant} Let $p(x)$ and $q(x)$ be two
polynomials over an algebraicly closed field. We define 
{\em generalized resultant} of $p(x)$ and $q(x)$ to be
$$
R(\lambda,\mu)=\prod_{i,j}{(\lambda \alpha_i+\mu \beta_j)}
$$
where $\left\{\alpha_i\right\}$ and $\left\{\beta_j\right\}$ are roots of
polynomials $p(x)$ and $q(x)$ respectively.
\end{definition}

Generalized resultant is a polynomial in two variables. It is easy to show that
its coefficients are polynomials in coefficients of $p(x)$ and $q(x)$ so the
condition on the base field to be algebraicly closed can be omitted.

\begin{theorem}The characteristic polynomial $\chi(\lambda,\mu,F)$ for algebra
$\Mat_n$ in point $F\in \Mat_n^*$ over the entire algebra $\Mat_n$ in basis of
matrix units is equal to the generalized resultant of 
characteristic polynomial of $F$ (as a matrix) with itself times $(-1)^{\frac{n(n-1)}{2}}$. That is
$$
\chi(\lambda,\mu,F)=(-1)^{\frac{n(n-1)}{2}}\det(F)(\lambda+\mu)^n\prod_{i\neq j}{(\lambda \alpha_i+\mu \alpha_j)}
$$
where $\alpha_i$ are eigenvalues of $F$ (this formulation assumes
that  the base field is algebraicly closed).
\end{theorem}
\begin{proof}
We will make use of theorem \ref{chi_invariance}. The coadjoint
action on $\Mat_n^*$ is simply conjugation by invertible matrices. The generic
orbit consists of diagonalizable matrices. Thus we can compute $\chi(\lambda,\mu)$
by assuming first that $F$ is diagonal and then extrapolating the resulting polynomial
to the case of all $F$.

Assume the base field to be $\mathbb C$.
Let $F=\diag(\alpha_1,...,\alpha_n)$. We choose a basis $\left\{E_{i.j}\right\}$
of matrix units in the algebra $\Mat_n$. The only case when $F(E_{i,j}E_{k,l})$
is non-zero is when $i=l$ and $j=k$. Thus the multiplication table $A$ (restricted
to the subspace of diagonal matrices in $\Mat_n^*$) is
$$
\begin{array}{ccccccccccccccccc}
& \vline &  & E_{i,i}   &  &\vline& & E_{i,j}^+ & &\vline & & E_{i,j}^- & \cr
\hline
& \vline& \alpha_1  & & 0& \vline &&&& \vline \cr
E_{i,i} & \vline&  & \ddots  & &\vline& & 0 & &\vline & & 0\cr
& \vline& 0&  & \alpha_n &\vline& & & & \vline & \cr
\hline
& \vline& & && \vline & & & & \vline &\alpha_{j'} & & 0 \cr
E_{i,j}^+ & \vline & & 0 & &\vline & & 0 & &\vline &&\ddots	\cr
& \vline& & && \vline & & & & \vline &0 & &\alpha_{j''}\cr
\hline
& \vline& & && \vline &\alpha_{i'} & &0 & \vline & \cr
E_{i,j}^- &\vline & & 0 & &\vline & &\ddots & &\vline & & 0\cr
& \vline& & && \vline &0 & &\alpha_{i''} & \vline & \cr
\end{array}
$$
here $E_{i,j}^+$ denotes elements $E_{i,j}$ with $i>j$ and $E_{i,j}^-$ denotes
elements $E_{i,j}$ with $i<j$. The matrix $\lambda A+\mu A^T$ will have 
$(\lambda+\mu)\alpha_i$ in the $E_{i,i}\times E_{i,i}$ block, and the pair 
$(E_{i,j}^+,E_{i,j}^-)$ will produce a $2\times 2$ matrix
$$
\left(
\begin{array}{cc}
0 & \lambda\alpha_j+\mu\alpha_i \cr
\lambda\alpha_i+\mu\alpha_j & 0 \cr
\end{array}
\right)
$$
Computing the determinant yields
$$
(-1)^{\frac{n(n-1)}{2}}(\lambda+\mu)^n\prod_i\alpha_i\prod_{i\neq j}{(\lambda\alpha_i+\mu\alpha_j)}
=(-1)^{\frac{n(n-1)}{2}}\prod_{i,j}{(\lambda\alpha_i+\mu\alpha_j)}
$$
thus proving the theorem for the case when $F$ is diagonal. But characteristic
polynomial $\det(F-x)$ is invariant under coadjoint action. Thus this expression
is true for all $F$ up to a possibly missing factor depending only on $F$ (but not
$\lambda$ or $\mu$) which is quasi-invariant under co\-adjoint action. However,
in view of the fact that
this multiple must be a polynomial in $F$ and that the degree of the 
expression above in $F$ is exactly $n^2$ this multiple must be trivial.

The case of an arbitrary field is proved by observing that both sides of the equality
are polynomials with integral coefficients and thus if equality holds over $\mathbb C$
it should hold over any field.
\end{proof}
One easily observes that for any $F$ with all distinct, non-zero eigenvalues (as a matrix)
the characteristic polynomial does not vanish.
\end{example}
\subsection{Type 2}
The easiest example of a type 2 algebra is given by direct sum of a type 1 algebra
with an algebra with trivial (identically $0$) multiplication law.

A more interesting example is given by the following construction:

\begin{example}\label{abc=0}
Let $V$ be a vector space of dimension $k$ and let $B:V\times V\rightarrow W$ be a bilinear
map of crossproduct $V\times V$ into vector space $W$ of dimension $m$.

We define algebra $\mathfrak A(B)$ by the following multiplication table:

$$
\begin{array}{c|cc}
 & V & W \\
 \hline 
 V & B & 0 \\
 W & 0 & 0 \\
\end{array}
$$

The algebra $\mathfrak A(B)$ possesses a remarkable property - the product of any three elements
is always zero. Thus it is always associative, no matter what $B$ is.

It is straightforward to see that for any $F\in \mathfrak A(B)^*$ we have $W\subset \Nil_F$.
If one chooses $B$ and $F\in\mathfrak A(B)^*$ in such a way that $\det(F(B))$ is non-zero we 
obtain an example of a type 2 pair $(\mathfrak A(B), F)$. 

Since the inequality $\det(F(B))\neq 0$ defines a Zariski open set of functionals $F$ any
algebra $\mathfrak A(B)$ that possesses $F$ of type 2 is a type 2 algebra.
\end{example}

\subsection{Type 3}

Let us consider a special case of the example \ref{abc=0} with $\dim W=1$.

In this situation the matrix $\AF$ depends only on the value of $F$ on the 
single basis vector of $W$ and $B$. Let us select $0\neq w\in W$ and any $F$ such that
$F(w)=1$.

By manipulating $B$ we can thus set $\AF$ to anything we like with only restriction that the
last row and column are identically $0$. This provides a lot of examples of type 3 algebras, 
in particular the following $B$ will do just fine:
$$
B=\left(
\begin{array}{ccc}
0 & 0 & 1 \\
0 & 0 & 1 \\
0 & 0 & 0 \\
\end{array}
\right)
$$

For this $B$ the space $\Stab_F(0)=\ker \AF$ is spanned by $v_1, v_2 \textrm{ and } w$ and
$\Stab_F(\infty)=\ker \AF^T$ is spanned by $v_1-v_2, v_3 \textrm{ and } w$. Thus $\Nil_F=W$
but $\det\left(\lambda B+\mu B^T\right)=0$.

\section{Spaces $\Stab_F(\alpha)$}
\begin{definition}
Let $\mathfrak A$ be an associative algebra and $F$ be a linear functional on it. 
We define
$$
\Stab_F(\alpha):=\left\{a\in\mathfrak A: \forall x\in\mathfrak A\follows F(ax)-\alpha F(xa)=0\right\}
$$
In other words 
$$
\Stab_F(\alpha)=\ker \left(\AF-\alpha \AF^T\right)
$$
\end{definition}

If one considers Lie algebra $\mathfrak A$ with bracket $[a,b]=ab-ba$ then $\Stab_F(1)=\Stab_F$
in the conventional definition of stabilizer of a linear functional on a Lie algebra.

\begin{example} Returning to the example \ref{Mat_n} we see that for $i\neq j$
$$
\Stab_F\left(\frac{\alpha_i}{\alpha_j}\right)=\mathbb C \cdot e_{ij}
$$
and 
$$
\Stab_F(1)=\span<e_{11},...,e_{nn}>
$$
\end{example}

\begin{theorem}\label{stab_mult}
$$
\begin{array}{rcl}
\Stab_F(\alpha) \cdot \Stab_F(\beta)&\subset& \Stab_F(\alpha\beta) \\
\Stab_F(0)\cdot \Stab_F(\infty)&\subset& \Nil_F \\
\Stab_F(0)\cdot \mathfrak A &\subset& \Stab_F(0) \\
\mathfrak A \cdot \Stab_F(\infty) &\subset& \Stab_F(\infty)\\
\dim \Stab_F(\alpha)&=&\dim \Stab_F(1/\alpha)
\end{array}
$$
\end{theorem}
\begin{proof}
Let $a\in\Stab_F(\alpha)$ and $b\in\Stab_F(\beta)$. Then for all $x$
$$
F((ab)x)=F(abx)=\alpha F(bxa)=\alpha\beta F(x(ab))
$$
For $a\in\Stab_F(\infty)$ and $b\in\Stab_F(\beta)$, $\beta\neq 0$ and any $x$ we have:
$$
F(x(ab))=F(xab)=\frac{1}{\beta}F(bxa)=0
$$
and
$$
F(x(ba))=F(xba)=F((xb)a)=0
$$
Secondly, for $a\in\Stab_F(0)$ and $b\in\Stab_F(\infty)$ and any $x$ we have
$$
F((ab)x)=F(a(bx))=0
$$
and
$$
F(x(ab))=F((xa)b)=0
$$
Thus $\Stab_F(0)$ is a right ideal in $\mathfrak A$ and $\Stab_F(\infty)$ is a 
left ideal in $\mathfrak A$. The product of $\Stab_F(0)$ and $\Stab_F(\infty)$
must be in $\Stab_F(0)\intersect\Stab_F(\infty)=\Nil_F$.

Lastly from linear algebra we know that for any matrix $R$ 
$$
\dim \ker R=\dim \ker R^T
$$
Thus, from definition of $\Stab_F(\alpha)$, it follows that $\dim \Stab_F(\alpha)=\dim \Stab_F(1/\alpha)$.
\end{proof}
{\bf Corollary:} $\Stab_F(1)$, $\Stab_F(0)$ and $\Stab_F(\infty)$ are subalgebras of $\mathfrak A$.

\begin{theorem}Let $\mathfrak A$ be a unital associative algebra and $F$ a linear
functional on it. Then for all
$\alpha\neq 1$ we have
$$
F(\Stab_F(\alpha))=\left\{0\right\}
$$
\end{theorem}
\begin{proof}
Indeed, consider first the case when $\alpha$ is finite.
By definition for any element $a\in\Stab_F(\alpha)$ we have
$$
0=F(a\cdot 1)-\alpha F(1\cdot a)=(1-\alpha)F(a)
$$
and thus $F(a)=0$.

Similarly, for $\alpha=\infty$ we must have
$$
0=F(1\cdot a)=F(a)
$$
\end{proof}

\section{Multiplicative functionals}
Multiplicative functionals play an important role in classical representation theory
and in the study of commutative algebras.
\begin{proposition} Let $F$ be a multiplicative functional on unital associative algebra $\mathfrak A$. 
Then $\AF$ has rank $1$.
\end{proposition}
\begin{proof}
Indeed, from definition
$$
\AF=\left(F(e_ie_j)\right)_{i,j=1}^n=\left(F(e_i)F(e_j)\right)_{i,j=1}^n=\left(F(e_i)\right)_{i=1}^n\left(F(e_j)\right)_{j=1}^n
$$
Since $F(1)=1$ the matrix $\AF$ cannot be $0$.
\end{proof}

\begin{lemma}
Let $\mathfrak A$ be an associative algebra with linear functional $F$ such that
$\Stab_F(0)=\Stab_F(\infty)=\Nil_F$. Then $\Nil_F$ is an ideal.
\end{lemma}
\begin{proof}
This is a direct consequence of the corollary to theorem \ref{stab_mult}.
\end{proof}

\begin{theorem} Let $F$ be a linear functional on unital associative algebra $\mathfrak A$
such that the matrix $\AF$ has rank $1$ and $F(1)=1$.
Then $F$ is multiplicative.
\end{theorem}
\begin{proof}
Since rank $\AF$ is $1$ it must be that $\codim \Stab_F(0)=\codim \Stab_F(\infty)=1$.
On the other hand 
$$
F(x)=F(1\cdot x)=F(x \cdot 1)
$$ 
and thus 
$$
F(\Stab_F(0))=F(\Stab_F(\infty))=0
$$
But $\codim \ker F = 1$. 

Thus $\Stab_F(0)=\Stab_F(\infty)=\ker F=\Nil_F$. From the previous lemma we know that
$\ker F=\Nil_F$ is an ideal.

Thus for any $a,b\in\mathfrak A$:
$$
\begin{array}{rcl}
F(ab)&=&F((F(a)+(a-F(a)))(F(b)+(b-F(b))))=\\
&=&F(F(a)F(b))+F(F(a)(b-F(b)))+F((a-F(a))F(b))+\\
&&+F((a-F(a))(b-F(b)))\\
\end{array}
$$
The first term is exactly $F(a)F(b)$. The second and third terms vanish because
$x-F(x)$ belongs to $\ker F$. The last term also vanishes because $\ker F$ is an ideal.
Thus $F(ab)=F(a)F(b)$ and $F$ is multiplicative.
\end{proof}

It is interesting to note that this proof applies to both finite- and infinite-dimensional
associative algebras. If the algebra has been endowed with topology than we 
can restrict out attention to only continuous functionals.

\section{Regular functionals}
The multiplicative functionals correspond to the situation when $\AF$ (or $\lambda \AF+\mu \AF^T$)
has the smallest rank possible. The case when $\AF$ (or $\lambda \AF+\mu \AF^T$)
has the maximum possible rank is described by {\em regular} functionals.

Following \cite{D1} we first prove the following lemma:
\begin{lemma}\label{Transversality lemma for constant lambda and mu}
Let $W$ be a vector subspace of $\mathfrak A$. Fix $\lambda$ and $\mu$. Then, the set $R$ of all $F\in \mathfrak A^*$ such
that $W \intersect \ker \left(\lambda \AF+\mu \AF^T\right) \neq 0$ is closed in
$\mathfrak A^*$.
\end{lemma}
\begin{proof}
Let $e_1,\ldots,e_n$ be the basis of $\mathfrak A$ such that the first $p$ vectors
form the basis of $W$. The system of equations
$$
\left(\lambda \AF + \mu \AF^T\right)\sum_{i=1}^p \epsilon_i e_i=0
$$
is equivalent to
$$
\sum_{i=1}^{p}\epsilon_i\left(\lambda F(e_ke_i)+\mu F(e_ie_k)\right)=0
$$
where $k$ ranges from $1$ to $n$. The matrix of the latter system has entries that 
are linear in $F$. The existence of the non-zero solution is equivalent to the requirement
that all $p$-minors vanish. This proves that the set $R$ can be defined as a solution to
the system of polynomial equations. Hence it is closed.
\end{proof}

This argument can also be generalized to the case of $\lambda$ and $\mu$ varying with $F$. 
However, one would have to specify the degree of regularity of $\lambda$ and $\mu$.
We prefer to state it with the assumption of continuity:

\begin{lemma}\label{Transversality lemma for continuous lambda and mu}
Let $W$ be a vector subspace of $\mathfrak A$. Fix two continuous functions $\lambda(F)$ and $\mu(F)$. The set $R$ of all $F\in \mathfrak A^*$ such
that 
$$
W \intersect \ker \left(\lambda(F) \AF+\mu(F) \AF^T\right) \neq \left\{0\right\}
$$
 is closed in $\mathfrak A^*$.
\end{lemma}

\begin{theorem}
Let $S$ be a subspace of $\mathfrak A^*$. 
Let $F$ be such that for a fixed pair of continuously differentiable functions $\lambda(F)$, $\mu(F)$ the dimension of
 the space 
 $$\ker\left(\lambda(F) \AF+\mu(F)\AF^T\right)$$ 
 is the smallest among functionals
in small neighbourhood of $F$ inside affine set $F+S$.

Then for any $G\in S$, any $x\in \ker\left(\lambda(F)\AF+\mu(F)\AF^T\right)$ and 
any $y\in \ker\left(\mu(F)\AF+\lambda(F)\AF^T\right)$
we have
$$
G\left(\lambda(F)xy+\mu(F)yx\right)+\left(L_G\lambda\right)(F)F(xy)+\left(L_G\mu\right)(F)F(yx)=0
$$
Here $L_G$ denotes directional derivative.
\end{theorem}
\begin{proof}
Pick any $G\in S$. 
Let $F_\epsilon=F+\epsilon G$. Choose any subspace $W\subset\mathfrak A$ of complementary
dimension that is
transversal to $\ker\left(\lambda(F)\AF+\mu(F)\AF^T\right)$. The set of all $\epsilon$
for which $V_\epsilon=\ker\left(\lambda(F_{\epsilon})\AFVA{F_{\epsilon}}+\mu(F_{\epsilon})\AFVA{F_\epsilon}^T\right)$
remains transversal to $W$ is open by lemma \ref{Transversality lemma for constant lambda and mu} and thus
contains a small neighbourhood of $0$. Let $\left\{e_i\right\}$ be the basis of $\mathfrak A$
with first $p$ vectors forming the basis of $W$.

Pick $x\in  \ker\left(\lambda(F)\AF+\mu(F)\AF^T\right)$. As $W$ and $\ker\left(\lambda(F_{\epsilon})\AFVA{F_{\epsilon}}+\mu(F_{\epsilon})\AFVA{F_\epsilon}^T\right)$
form a decomposition of $\mathfrak A$ there exists $w\in W$ such that
$$
x-w \in \ker\left(\lambda(F_{\epsilon})\AFVA{F_{\epsilon}}+\mu(F_{\epsilon})\AFVA{F_\epsilon}^T\right)
$$
I.e. for all $i$ 
$$
\lambda(F_{\epsilon})F_{\epsilon}((x-w)e_i)+\mu(F_{\epsilon})F_{\epsilon}(e_i(x-w))=0
$$
Hence
$$
\lambda(F_{\epsilon})F_{\epsilon}(we_i)+\mu(F_{\epsilon})F_{\epsilon}(e_iw)=\lambda(F_{\epsilon})F_{\epsilon}(xe_i)+\mu(F_{\epsilon})F_{\epsilon}(e_ix)
$$
The homogeneous part of this system is the linear system that defines $W\intersect \ker\left(\lambda(F_{\epsilon})\AFVA{F_{\epsilon}}+\mu(F_{\epsilon})\AFVA{F_\epsilon}^T\right)$.
As this intersection is trivial the rank of it must be $p$ for all $\epsilon$ in a small neighbourhood of $0$. Noticing that the right hand
side can be derived by substituting $w=x$ (i.e. our system is of the type $Qw=Qx$) we conclude
that there is a unique solution $w_{\epsilon}$. Recall that $F_{\epsilon}$ is linear in $\epsilon$. 
When $\lambda(F)$ and $\mu(F)$ are constant $w_\epsilon$ is a rational function in $\epsilon$.
When $\lambda(F)$ and $\mu(F)$ are continuously differentiable the non-vanishing $p$-minors of the homogeneous part of the system 
remain non-vanishing for $\epsilon$ in a small neighbourhood of $0$ and thus $w_{\epsilon}$ is continuously differentiable as well.

Pick now any $y\in\ker\left(\lambda(F)\AF+\mu(F)\AF^T\right)$. We have:
$$
\lambda(F_{\epsilon})F_{\epsilon}((x-w_{\epsilon})y)+\mu(F_{\epsilon})F_{\epsilon}(y(x-w_{\epsilon}))=0
$$
Differentiating with respect to $\epsilon$ produces:
$$
\begin{array}{c}
\lambda(F_{\epsilon})G_{\epsilon}((x-w_{\epsilon})y)+\mu(F_{\epsilon})G_{\epsilon}(y(x-w_{\epsilon}))-
\left(\lambda(F_{\epsilon})F_{\epsilon}(w_{\epsilon}'y)+\mu(F_{\epsilon})F_{\epsilon}(yw_{\epsilon}')\right)+\\
+\lambda(F_{\epsilon})'F_{\epsilon}((x-w_{\epsilon})y)+\mu(F_{\epsilon})'F_{\epsilon}(y(x-w_{\epsilon}))=0
\end{array}
$$
After setting $\epsilon=0$ we have:
$$
\lambda(F)G(xy)+\mu(F)G(yx)+\left.\lambda(F_{\epsilon})'\right|_{\epsilon=0}F(xy)+\left.\mu(F_{\epsilon})'\right|_{\epsilon=0}F(yx)=0
$$
and using notation of Lie derivative we get:
$$
\lambda(F)G(xy)+\mu(F)G(yx)+\left(L_G\lambda\right)(F)F(xy)+\left(L_G\mu\right)(F)F(yx)=0
$$
\end{proof}

{\bf Corollary 1.} For the case of constant $\lambda$ and $\mu$ we obtain a simpler expression:
$$
\lambda(F)G(xy)+\mu(F)G(yx)=0
$$

{\bf Corollary 2.} Let $F$ be a functional such that $\Stab_F(\alpha(F))$ has the smallest
dimension in a neighbourhood of $F$. ($\alpha(\cdot)$ is fixed, finite and non-zero in $F$). Then for all $x\in \Stab_F(\alpha(F))$
and $y\in\Stab_F(1/\alpha(F))$ we have:
$$
xy-\alpha yx=0
$$

{\bf Corollary 3.} Let $F$ be a functional with the smallest dimension of $\Stab_F(1)$.
Then $\Stab_F(1)$ is a commutative subalgebra of $\mathfrak A$.

{\bf Corollary 4.} Let $F$ be a functional with the smallest dimension of $\Stab_F(0)$.
Then 
$$
\Stab_F(0)\cdot \Stab_F(\infty)=\left\{0\right\}
$$
In this case $\Nil_F$ is a subalgebra of $\mathfrak A$ with trivial (identically $0$)
multiplication law.

{\bf Corollary 5.} Let $F$ be a functional such that $\Stab_F(\alpha(F))$ has the smallest
dimension in a neighbourhood of $F$. ($\alpha(\cdot)$ is fixed, finite and non-zero). Then the set
$$
\begin{array}{l}
\left[\Stab_F(\alpha(F)),\Stab_F(1/\alpha(F))\right]_{\alpha(F)}:=\\
\qquad\qquad:=\left\{xy-\alpha(F) yx:x\in\Stab_F(\alpha(F))\textrm{ and } y\in\Stab_F(1/\alpha(F))\right\}
\end{array}
$$
has at most dimension $1$.
\begin{proof}
Let $\lambda=1$ and $\mu=-\alpha(F)$. From the theorem we know that
$$
\lambda G(xy)+\mu G(yx)+\left(L_G\lambda\right)(F)F(xy)+\left(L_G\mu\right)(F)F(yx)=0
$$
Substituting we get
$$
G(xy)-\alpha(F)G(yx)-\left(L_G\alpha\right)(F)F(yx)=0
$$
There exists an element $z\in\mathfrak A$ such that for all $G\in\mathfrak A^*$ we have 
$$
\left(L_G\alpha\right)(F)=G(z)
$$
Thus for all $G\in\mathfrak A^*$
$$
G(xy)-\alpha(F)G(yx)-G(z)F(yx)=0
$$
and
$$
G\left(xy-\alpha(F)yx-zF(yx)\right)=0
$$
which implies
$$
xy-\alpha(F)yx-zF(yx)=0
$$
and thus $\dim \left[\Stab_F(\alpha(F)),\Stab_F(1/\alpha(F))\right]_{\alpha(F)}$ is at most $1$.
\end{proof}
\begin{definition} Let $\mathfrak A$ be an associative algebra and $F$ a linear functional.
Let $S$ be a subspace of linear functionals on $\mathfrak A$. Let $\alpha$ be either a constant
or a continuously differentiable function. 

We will call $F$ $(\alpha,S)$-regular if the space $\Stab_F(\alpha(F))$ has the smallest
dimension among functionals in the neighbourhood of the affine set $F+S$.

In the case $S=\mathfrak A^*$ we will simply call $F$ $\alpha$-regular.
\end{definition}

\begin{definition} Let $\mathfrak A$ be an associative algebra and $F$ be a $1$-regular
linear functional on it.

We define $\ind \mathfrak A:=\dim \Stab_F(1)$
\end{definition}

{\bf Note:} the index defined above is the same as the index of Lie algebra $\mathfrak A^{Lie}$
obtained from $\mathfrak A$ by defining the bracket operation as $[a,b]=ab-ba$.

\section{Type 1 algebras}
We will now study type 1 algebras in more detail. 

We will establish a criteria for recognizing type $1$ algebras, analyze their characteristic polynomial, 
obtain decomposition into subspaces $V(\alpha)$ and study tensor products of type $1$ algebras.

In the end we will prove the identity
$$
\ind \Mat_n \otimes \mathfrak B=n
$$
for any finite dimensional unital associative algebra $\mathfrak B$ over complex numbers that has index $1$.
\subsection{Recognizing type 1 algebras}
\begin{definition}
For an associative algebra $\mathfrak A$ and a functional $F$ we define the skew-symmetric
form $B_F$ as
$$
B_F(a,b)=F(ab-ba)
$$
and the symmetric form $Q_F$ as
$$
Q_F(a,b)=F(ab+ba)
$$
\end{definition}

\begin{theorem}\label{q_non_degenerate}Let $\mathfrak A$ be an associative algebra and $F$ a linear functional
on it. Suppose that
the restriction of the form $Q_F$ on the space $\Stab_F(1)$ is a non-degenerate
symmetric form. Then $F$ is of type 1.
\end{theorem}
\begin{proof}

By definition $F$ is of type 1 if and only if $\det(\lambda\AF+\mu\AF^T)\neq 0$.
Let $\epsilon=(\lambda+\mu)/2$ and $\sigma=(\lambda-\mu)/2$. Let $k$ denote $\dim\Stab_F(1)$.
We compute:
$$
\begin{array}{rcl}
\det(\lambda\AF+\mu\AF^T)&=&\det((\epsilon+\sigma)\AF+(\epsilon-\sigma)\AF^T)=\\
\qquad&=&\det(\epsilon(\AF+\AF^T)+\sigma(\AF-\AF^T))=\\
\qquad&=&\det(\epsilon Q_F+\sigma B_F)=\\
\qquad&=&\epsilon^{k}\det\left(\left.Q_F\right|_{\Stab_F(1)}\right)\det\left(\sigma \left.B_F\right|_{\Stab_F(1)^{\perp}}\right)+o(\epsilon^{k})
\end{array}
$$
Since $B_F$ is non-degenerate on $\Stab_F(1)$-transversal subspace of $\mathfrak A$ of 
complimentary dimension we must have $det(\lambda\AF+\mu\AF^T)\neq 0$.
\end{proof}

The following is an easy consequence of theorem \ref{q_non_degenerate}:
\begin{theorem} Let $\mathfrak A$ be a unital associative algebra of index $1$.
(i.e. $\ind \mathfrak A=1$). Then $\mathfrak A$ is type 1.
\end{theorem}
\begin{proof}
Indeed, consider the set of $1$-regular functionals $F$ that do not vanish on unity.
The form $Q_F$ is scalar and equal to $F(1)$. Thus all such functionals $F$ are type 1.
\end{proof}

We are now in position to formulate sufficient condition for an algebra to be type 1:
\begin{theorem} \label{F is 1-regular}
Let $\mathfrak A$ be an associative algebra and let $F$ be a linear functional on it
such that $\Stab_F(1)$ is commutative and the restriction of $Q_F$ on $\Stab_F(1)$
is non-degenerate. Then $F$ is $1$-regular and $\mathfrak A$ is a type 1 algebra.
\end{theorem}
\begin{proof}
If $F$ is $1$-regular then the fact that $\mathfrak A$ is type 1 follows from theorem
\ref{q_non_degenerate}. 

We now concentrate on proving that $F$ is $1$-regular.

Recall the following facts from the theory of Lie algebras:

The map $\mathfrak A \rightarrow \mathfrak A^*$ defined as
$$
\ad{}^*(a):F(x)\mapsto F(ax-xa)
$$
vanishes exactly on $\Stab_F(1)$. The images in each functional $F$ form a locally 
integrable distribution. 

Under the map $\ad{}^*$ the form $B_F$ is mapped to Kirillov's form on the image
of $\ad{}^*$. 

Kirillov's form provides symplectic structure on the leaves of this distribution.
Of particular note to us is the fact that for any two functionals $F$ from the same
leaf the dimension of $\Stab_F(1)$ is the same.

Thus if we were to construct a section in $F$ - that is a manifold of dimension
complimentary to the dimension of the leaf that $F$ belongs to, containing $F$ and transversal
to the leaf containing $F$ - and prove that 
in a small neighbourhood in the section around $F$ the dimension of $\Stab_F(1)$ 
does not vary we would be able to parameterize all points in the neighbourhood of
$F$ by the leaf that point belongs to and the intersection of that leaf with 
the section in $F$. Since in a small neighbourhood of $F$ inside the section the
dimension of $\Stab_F(1)$ does not vary we would obtain that $F$ is $1$-regular.

Consider the following map of $\Stab_F(1)$ into $\mathfrak A^*$:
$$
\rho:a\mapsto \tilde{F}(x):=F(x)+F(xa)+F(ax)=F(x)+Q_F(a,x)
$$
This map possesses the following properties:
\begin{itemize}
\item the image of $\rho$ is a linear submanifold of $\mathfrak A^*$
\item map $\rho$ is a bijection between $\Stab_F(1)$ and $\image \rho$. Indeed,
assume this is not true and there are two $a$ and $b$ that mapped into the
same functional $\tilde{F}$. Then, for all $x$:
$$
F(x)+Q_F(a,x)=F(x)+Q_F(b,x)
$$
Thus
$$
Q_F(a-b,x)=0
$$
Since $Q_F$ is non-degenerate on $\Stab_F(1)$ we must have $a-b=0$ and thus $a$ and $b$
coincide.
\item the tangent space of $\image \rho$ in point $F$ is transversal to the tangent
space of the leaf passing through $F$. Indeed, assume this is not so and there exists
$a\in\Stab_F(1)$ and $b\in\mathfrak A$ such that
$$
Q_F(a,x)=F(bx-xb)
$$
We note that $a$ cannot be zero as otherwise the tangent vector itself must be $0$.
Since $Q_F$ is non-degenerate there exists $c\in\Stab_F(1)$ such that $Q_F(a,c)\neq 0$.
By definition of $\Stab_F(1)$ we have $F(bc-cb)=0$ - a contradiction.
Therefore these spaces are transversal.
\item for all $\tilde{F}$ we have $\Stab_F(1)\subset \Stab_{\tilde{F}}(1)$. Indeed,
let $x$ be an arbitrary element of $\mathfrak A$, $b$ an arbitrary element of $\Stab_F(1)$
and $a$ be an element of $\Stab_F(1)$ that defines $\tilde{F}$. 

By the assumption of the theorem $\Stab_F(1)$ is commutative, thus
$$
\tilde{F}(bx)=F(bx)+F(bxa)+F(abx)=F(xb)+F(xba)+F(axb)=\tilde{F}(xb)
$$
\item there exists a small neighbourhood of $F$ inside $\image \rho$ where the 
dimension of $\Stab_F(1)$ does not vary. Indeed, the set of all functionals $G\in\image\rho$ such
that $\dim \Stab_G(1)$ is greater than $\dim \Stab_F(1)$ is defined by a set of
polynomial equations in $G$ (namely that the minors of the order $\dim \mathfrak A-\dim \Stab_F(1)$ of the matrix $G(e_ie_j-e_je_i)$
all vanish) and thus is Zariski closed. The complement is a Zariski open set - and
$F$ already belongs to it. Since we have already proven that all $\tilde{F}$ in $\image \rho$
have $\dim\Stab_{\tilde{F}}(1)\ge \dim\Stab_F(1)$, the set of all $\tilde{F}$ such 
that $\dim\Stab_{\tilde{F}}(1)=\dim\Stab_F(1)$ is Zariski open.
\end{itemize}
Therefore, $\image \rho$ is the section we desire and thus $F$ is $1$-regular.
\end{proof}

\subsection{Characteristic polynomial of type $1$ algebras}

\begin{theorem}
$\chi_F(\lambda,\mu)$ is divisible by $(\lambda-\mu/\alpha)^{\dim\Stab_F(\alpha)}$  ($\mu^{\dim\Stab_F(0)}$ for the
case $\alpha=0$).
\end{theorem}
\begin{proof}
Let $\epsilon=\lambda-\mu/\alpha$ and $\sigma=\mu$. Let $k$ denote $\dim\Stab_F(\alpha)$. We will choose two basises
in $\mathfrak A$: one such that the first $k$ vectors form $\Stab_F(\alpha)$ and second so that the first $k$ vectors
form $\Stab_F(1/\alpha)$.
We compute:
$$
\begin{array}{l}
\det(\lambda\AF+\mu\AF^T)=\det((\epsilon-\sigma/\alpha)\AF+(\sigma\AF^T)=\\
\qquad=\det(\epsilon\AF+\sigma(\AF-\alpha\AF^T)/\alpha)=\\
\qquad=\epsilon^{k}\det\left(\left.\AF\right|_{\Stab_F(\alpha)\times\Stab_F(\frac{1}{\alpha})}\right)\det\left(\left.\frac{\sigma(\AF-\alpha\AF^T)}{\alpha}\right|_{\Stab_F(\alpha)^{\perp}\times\Stab_F(\frac{1}{\alpha})^{\perp}}\right)+o(\epsilon^{k})
\end{array}
$$
Thus $\det(\lambda\AF+\mu\AF^T)$ is divisible by at least $\epsilon^{k}$.

The cases $\alpha=0$ and $\alpha=\infty$ are resolved in a similar manner:
$$
\begin{array}{l}
\det(\lambda\AF+\mu\AF^T)=\\
\qquad\qquad=\lambda^{\dim\ker\AF^T}\det\left(\left.\AF\right|_{\Stab_F(\infty)\times\Stab_F(0)}\right)\det\left(\left.\AF^T\right|_{\Stab_F(\infty)^{\perp}\times\Stab_F(0)^{\perp}}\right)+\\
\qquad\qquad\qquad+o(\epsilon^{\dim\ker\AF^T})=\\
\qquad\qquad=\mu^{\dim\ker\AF}\det\left(\left.\AF^T\right|_{\Stab_F(0)\times\Stab_F(\infty)}\right)\det\left(\left.\AF\right|_{\Stab_F(0)^{\perp}\times\Stab_F(\infty)^{\perp}}\right)+\\
\qquad\qquad\qquad+o(\epsilon^{\dim\ker\AF})\\
\end{array}
$$
\end{proof}

{\bf Corollary 1.} We see from the proof that $\dim\Stab_F(\alpha)$ coincides with the highest degree $k$ such that
$\chi_{F}(\lambda,\mu)$ is divisible by $\left(\lambda-\mu/\alpha\right)^k$ if and only if the restriction of
the form $\AF$ ($\AF^T$ for $\alpha=0$) on the space $\Stab_F(\alpha)\times\Stab_F(1/\alpha)$ is non-degenerate.

{\bf Corollary 2.} $\dim \Stab_F(1)$ is equal to the highest power of $\lambda-\mu$ that divides
$\chi_{F}(\lambda,\mu)$ if and only if the restriction of the form $Q_F$ to $\Stab_F(1)$ is
non-degenerate.

\begin{definition} We will call a functional $F$ on $\mathfrak A$ $\alpha$-precise if the
dimension of $\Stab_F(\alpha)$ is equal exactly to the highest power of $\left(\lambda-\mu/\alpha\right)$
that divides $\chi_{F}(\lambda,\mu)$.
\end{definition}

\begin{definition} Fix a continuous function $\alpha(F)$ defined on an open dense subset of $\mathfrak A^*$.
We call an associative algebra $\mathfrak A$ $\alpha(F)$-precise if there exists an open dense subset of functionals
$F$ that are $\alpha(F)$-precise.
\end{definition}

\begin{example} Let us compute characteristic polynomial for $(2,1;1,2)$ seaweed 
algebra\footnote{For definition of seaweed algebras see \cite{DK}} with 
multiplication table:
$$
\begin{array}{c|ccccc}
& a & b & c & d & e \\
\hline
a & a & b & 0 & 0 & 0 \\
b & 0 & 0 & b & 0 & 0 \\
c & 0 & 0 & c & 0 & 0 \\
d & 0 & 0 & d & 0 & 0 \\
e & 0 & 0 & 0 & d & e\\
\end{array}
$$
We compute
$$
\begin{array}{l}
\det\left(\lambda\AF+\mu\AF^T\right)=\lambda^2\mu^2b^2d^2(\lambda+\mu)(a+c+e)
\end{array}
$$
The space $\Stab_F(1)$ is generated by unity $a+c+e$, the space $\Stab_F(0)$
is spanned by $F(a)b-F(b)a$ and $F(e)d-F(d)e$ and the space $\Stab_F(\infty)$
is spanned by $F(b)c-F(c)b$ and $F(b)d-F(d)b$.

In this case we see that the dimensions of the spaces $\Stab_F(\alpha)$ match 
exactly the order of zero of polynomial $\chi_F(1,-x)$ in point $1/\alpha$.

Note that in this case $\mathfrak A=\Stab_F(0)\oplus\Stab_F(1)\oplus\Stab_F(\infty)$.
\end{example}
\subsection{Decomposition of type 1 algebras}
The previous example of $(2,1;1,2)$ seaweed algebra and the example of $n\times n$
matrices investigated in section \ref{mat_n} have the property that, for generic $F$,
the algebra $\mathfrak A$ is a direct sum of spaces $\Stab_F(\alpha)$.
While it is true that, for generic $F$, the spaces $\Stab_F(\alpha)$ do form a direct sum for any
type 1 algebra one can construct examples when their sum is not the entire algebra.

\begin{example}
The simplest way is to consider an algebra $\mathfrak A(B)$ from example~\ref{abc=0}
with $\dim W=1$ and extend it with unity:
$$
\begin{array}{c|ccc}
& 1 & V & w \\
\hline
1 & 1 & V & w \\
V & V & Bw & 0 \\
w & w & 0 & 0 \\
\end{array}
$$
The characteristic polynomial of such an algebra is equal to
$$
\chi_F(\lambda,\mu)=(\lambda+\mu)^2\det\left(\lambda B+\mu B^T\right) w^{\dim \mathfrak A}
$$
For this algebra $1$ and $w$ are always within $\Stab_F(1)$.

Without loss of generality we can assume that $F(1)=1$ and $F(w)\neq 0$. 
For $\alpha\neq 1$
the element $x=a\cdot 1 + v +b \cdot w$ belongs to $\Stab_F(\alpha)$ if and only if
$$
\left\{
\begin{array}{rcl}
a+F(v)+bF(w)&=&0\\
a(1-\alpha)F(v)+F(w)\left(B(v,\cdot)-\alpha B(\cdot, v)\right)&=&0\\
a(1-\alpha)F(w)&=&0\\
\end{array}
\right.
$$
Thus $a=0$ and $b=F(v)/F(w)$ and 
$$
(B-\alpha B^T)v=0
$$
Let us restrict out attention to those $B$ with $\det B\neq 0$. In this case the
previous equation reduces to $\left(1-\alpha B^{-1}B^T\right)v=0$ and the question
of whether $\mathfrak A$ is the direct sum of spaces $\Stab_F(\alpha)$ is equivalent
to the question of whether $B^{-1}B^T$ is diagonalizable.

Now consider the following family of matrices $B$:
$$
B=\left(\begin{array}{cc} 0 & B_0 \\ B_1 & 0 \end{array}\right)
$$
Then
$$
B^{-1}B^T=\left(\begin{array}{cc} 0 & B_1^{-1} \\ B_0^{-1} & 0 \end{array}\right)
\left(\begin{array}{cc} 0 & B_1^T \\ B_0^T & 0 \end{array}\right)=
\left(\begin{array}{cc} B_1^{-1}B_0^T & 0 \\ 0 & B_0^{-1}B_1^T \end{array}\right)
$$
There are plenty of choices for $B_0$ and $B_1$ that yield non-diagonalizable $B^{-1}B^T$.
\end{example}
However, one can extend the theory of Jordan decomposition of matrices to this case.
The following is a rather technical presentation of such.
\subsubsection{Definition of spaces $V_k(\alpha)$}
\begin{definition}
We define $V_{0}(\alpha)=\left\{0\right\}$ and $V_1(\alpha)=\Stab_F(\alpha)$.
\end{definition}

\begin{definition}\label{V_k} Let $\mathfrak A$ be an associative algebra. Fix $\alpha_0\neq\alpha$. 
We define $V_k(\alpha)$ - a space of "Jordan vectors" - as
$$
\begin{array}{l}
V_{k+1}(\alpha):=\left\{b\in\mathfrak A:\exists a\in V_k(\alpha) \imply\right.\\
\qquad\qquad \left. \imply\forall 
		x\in\mathfrak A\imply F(bx)-\alpha F(xb)=F(ax)-\alpha_0 F(xa)\right\}
\end{array}
$$
or in terms of multiplication table of $\mathfrak A$:
$$
V_{k+1}(\alpha):=\left\{b\in\mathfrak A:\exists a\in V_k(\alpha) \imply
    \AF b-\alpha\AF^T b= \AF a -\alpha_0 \AF^T a\right\}
$$
For $\alpha=\infty$ we define
$$
V_{k+1}(\infty):=\left\{b\in\mathfrak A:\exists a\in V_k(\alpha) \imply\left(\forall 
		x\in\mathfrak A\imply F(xb)=F(ax)-\alpha_0 F(xa)\right)\right\}
$$
or
$$
V_{k+1}(\infty):=\left\{b\in\mathfrak A:\exists a\in V_k(\alpha) \imply
    \AF^T b= \AF a-\alpha_0 \AF^T a\right\}
$$
\end{definition}

\begin{lemma} The spaces $V_k(\alpha)$ do not depend on the choice of $\alpha_0$.
\end{lemma}
\begin{proof}
Indeed this is so by definition for $k\leq 1$. Assume that for this statement holds
for all $k\leq n$.

Let $b$ be an element of $V_{n+1}(\alpha)$ constructed using $\alpha_0$. By definition, this
implies existence of $a \in V_n(\alpha)$ such that
$$
F(bx)-\alpha F(xb)=F(ax)-\alpha_0 F(xa)
$$
We will show existence of element $a''\in V_{n}(\alpha)$ such that $F(ax)-\alpha_0 F(xa)=F(a''x)-\alpha_1 F(xa'')$.

$$
\begin{array}{l}
F(ax)-\alpha_0 F(xa)=\\
\quad\quad\quad=(1-\delta)F(ax)+\delta\alpha F(xa) + \delta (F(a'x)-\alpha_1F(xa'))
    -\alpha_0 F(xa)=\\
\quad\quad\quad=(1-\delta)F(ax)+(\delta\alpha-\alpha_0)F(xa)+\delta(F(a'x)-\alpha_1 F(xa'))=\\
\quad\quad\quad=F(((1-\delta)a+\delta a')x)-\alpha_1 F(x ( (\alpha_0-\delta\alpha)a/\alpha_1 +\delta a'))
\end{array}
$$
Here $a'\in V_{n-2}(\alpha)$ satisfies $F(ax)-\alpha F(xa)=F(a'x)-\alpha_1 F(xa')$ by assumption of induction.

Choosing $\delta=\frac{\alpha_0-\alpha_1}{\alpha-\alpha_1}$ and $a''=(1-\delta)a+\delta a'$
we obtain $F(ax)-\alpha_1 F(xa)=F(a''x)-\alpha_0 F(xa'')$, where $a''$ is also in $V_{n-1}(\alpha)$.

Thus the space $V_{n+1}(\alpha)$ constructed using $\alpha_0$ is a subset of space $V_{n+1}(\alpha)$
constructed using any $\alpha_1\neq \alpha$. Therefore for any $\alpha_0$ and $\alpha_1$, different from
$\alpha$, the spaces $V_{n+1}(\alpha)$ are identical.

A similar argument can be used to prove the case $\alpha=\infty$. However, there is another way. We can
observe that when we introduce a new "transposed" multiplication law $a*b:=ba$ the parameter $\alpha$
is transformed into its inverse, i.e. space $\Stab_F(\alpha)$ become $\Stab_F(1/\alpha)$ and spaces
$V_k(\alpha)$ become spaces $V_k(1/\alpha)$. Since we already proven the case $\alpha=0$ we must conclude
that the case $\alpha=\infty$ is true as well.
\end{proof}
{\bf Note.} We emphasize that this definition is valid for {\em any} associative algebra, not necessarily 
finite dimensional or type 1.

For the case of finite dimensional type 1 algebras there is an equivalent way of defining $V_k(\alpha)$ that exposes their nature as Jordan spaces of an operator.
Choose $\alpha_0$ so that $\chi_F(\alpha_0)=\det\left(\AF-\alpha_0\AF^T\right)\neq0$. From definition,
$$
\begin{array}{l}
V_{k+1}(\alpha):=\left\{b\in\mathfrak A:\exists a\in V_k(\alpha) \imply
    \AF b-\alpha \AF^T b= \AF a-\alpha_0 \AF^T a \right\}=\\
\qquad\qquad =\left\{b\in\mathfrak A:\exists a\in V_k(\alpha) \imply
    \left(\AF -\alpha_0 \AF^T \right)^{-1}\left(\AF-\alpha \AF^T \right) b= a \right\}=\\
\qquad\qquad    =\left(\left(\left(\AF -\alpha_0 \AF^T \right)^{-1}\left(\AF-\alpha \AF^T \right)\right)^{k+1}\right)^{-1}V_0(\alpha)=\\
\qquad\qquad    =\left(\left(1-(\alpha-\alpha_0)\left(\AF -\alpha_0 \AF^T \right)^{-1}\AF^T\right)^{k+1}\right)^{-1}V_0(\alpha)=\\
\qquad\qquad    =\left(\left(\left(\AF -\alpha_0 \AF^T \right)^{-1}\AF^T-\frac{1}{\alpha-\alpha_0}\right)^{k+1}\right)^{-1}V_0(\alpha)
\end{array}
$$
We observe that $V_k(\alpha)$ is exactly the $k$-th level Jordan space of operator $\left(\AF -\alpha_0 \AF^T \right)^{-1}\AF^T$
corresponding to eigenvalue $\frac{1}{\alpha-\alpha_0}$. 

Thus
\begin{theorem} Let $\mathfrak A$ be a type 1 algebra. Then 
$$
\mathfrak A=\bigoplus_\alpha \bigcup_k V_k(\alpha)
$$
\end{theorem}

{\bf Remark.} For type 2 algebras there is no $\alpha_0$ such that $\AF-\alpha_0\AF^T$ is invertible.
However, by considering $\mathfrak A/\Nil_F$ instead of $\mathfrak A$ we notice that the induced 
form $\left(\AF-\alpha_0\AF^T\right)_{\Nil_F}$ is non-degenerate for most $\alpha_0$.
Therefore, $\mathfrak A/\Nil_F=\oplus_{\alpha}\cup_k \left(V_k(\alpha)/Nil_F\right)$.
For type 3 algebras one can construct an example where the spaces $V_k(\alpha)$ are pairwise transversal
for different values of $\alpha$, but do not form a direct sum.

\subsubsection{Properties of spaces $V_k(\alpha)$}

\begin{theorem}\label{V_k_prop} The spaces $V_k(\alpha)$ possess the following properties:
\begin{enumerate}
\item $V_k(\alpha) \subset V_{k+1}(\alpha)$
\item For $\alpha,\beta \notin \left\{0,\infty\right\}$ we have $V_k(\alpha) \cdot V_m(\beta) \subset V_{k+m-1}(\alpha\beta)$
\item For $\alpha\neq 0$ we have $V_k(\alpha)\cdot V_m(\infty)\subset V_{k+m-1}(\infty)$
\item For $\alpha\neq 0$ we have $V_k(\infty)\cdot V_m(\alpha)\subset V_{k+m-1}(\infty)$
\item For $\alpha\neq \infty$ we have $V_k(\alpha)\cdot V_m(0)\subset V_{k+m-1}(0)$
\item For $\alpha\neq \infty$ we have $V_k(0)\cdot V_m(\alpha)\subset V_{k+m-1}(0)$
\end{enumerate}

\end{theorem}
\begin{proof}
Property 1 follows by induction from the fact that
$V_{-1}(\alpha)\subset V_0(\alpha)$.

To prove property 2 we make induction on the parameter $N=k+m$.
The base of induction follows immediately from properties of $\Stab_F(\alpha)$ 
(theorem \ref{stab_mult}).
Assume that the statement is true for all $k$ and $m$ such that $k+m<N+1$.
For a given $\alpha$ and $\beta$ we pick $\alpha_0=0$ as this value is different from
both $\alpha$ and $\beta$.  Let $b_1\in V_k(\alpha)$, $b_2 \in V_m(\beta)$, where
$k+m=N+1$. Let $a_1\in V_{k-1}(\alpha)$ be an element corresponding to $b_1$
according to definition \ref{V_k} and $a_2\in V_{m-1}(\beta)$ be the element corresponding to $b_2$.
Let $x$ be an arbitrary element of $\mathfrak A$. 
Then:
$$
\begin{array}{l}
F(b_1b_2x)-\alpha\beta F(xb_1b_2)=\alpha F(b_2xb_1)-\alpha\beta F(xb_1b_2)+F(a_1b_2x)=\\
\qquad\qquad\qquad=\alpha F(a_2xb_1)+F(a_1b_2x)=\\
\qquad\qquad\qquad=\left( F(b_1a_2x)+F(a_1a_2x)\right)+F(a_1b_2x)=\\
\qquad\qquad\qquad=F\left(\left(b_1a_2+a_1a_2+a_1b_2\right)x\right)\\
\end{array}
$$
Now by assumption of induction we have 
$$
b_1a_2+a_1a_2+a_1b_2 \in V_{k+m-2}(\alpha\beta)
$$
and thus $b_1b_2$ is an element of $V_{k+m-1}(\alpha\beta)$.

Property 3. We perform induction the same way as in proof of property~2. For the same reasons we choose $\alpha_0=0$.  Let $b_1\in V_k(\alpha)$, $b_2 \in V_m(\infty)$, where
$k+m=N+1$. Let $a_1\in V_{k-1}(\alpha)$ be an element corresponding to $b_1$
according to definition \ref{V_k} and $a_2\in V_{m-1}(\infty)$ be the element corresponding to $b_2$. 

We compute:
$$
\begin{array}{c}
F(xb_1b_2)=F(a_2xb_1)=\frac{1}{\alpha}\left(F(b_1a_2x)+F(a_1a_2x)\right)=F\left(\frac{b_1a_2+a_1a_2}{\alpha}x\right)\\
\end{array}
$$
By assumption of induction we have 
$$
\frac{b_1a_2+a_1a_2}{\alpha}\in V_{k+m-2}(\infty)
$$
and thus $b_1b_2$ is an element of $V_{k+m-1}(\infty)$.

Property 4 is proved almost identically to property 3. We will write down the computation of $F(xb_1b_2)$:
$$
\begin{array}{c}
F(xb_1b_2)=\frac{F(b_2xb_1)+F(a_2xb_1)}{\alpha}=\frac{F(a_1b_2x)+F(a_1a_2x)}{\alpha}=F\left(\frac{a_1b_2+a_1a_2}{\alpha}x\right)\\
\end{array}
$$

Properties 5 and 6 can be proven by similar computation (it might be useful to use $\alpha=\infty$), however
we will simply refer to the correspondence $V_k(\alpha)\leftrightarrow V_k(1/\alpha)$ that occurs when
one considers a transposed algebra $\mathfrak A'$ with multiplication $a *b:=a\cdot b$.

\end{proof}

\begin{lemma}\label{kernel F}Let $\mathfrak A$ be an associative algebra with unity. Then
$F(V_k(\alpha))=0$ for all $\alpha \neq 1$.
\end{lemma}
\begin{proof}
Case $k=1$, $\alpha\neq \infty$: for all $b\in V_0(\alpha)$ we must have
$$
F(bx)=\alpha F(xb)
$$
Setting $x=1$ we get $F(b)=\alpha F(b)$, hence $F(b)=0$.

Case $k=1$, $\alpha=\infty$: we have $F(xb)=0$. Again setting $x=1$ yields
$F(b)=0$.

For arbitrary $k$ and $\alpha\neq \infty$ we proceed by induction. Again let us set $x=1$ in the definition
\ref{V_k}. We get
$$
F(b)-\alpha F(b)=F(a)
$$
But we already know that $F(a)=0$, thus $F(b)=0$ as well.

For arbitrary $k$ and  $\alpha=\infty$:  from the definition we get
$$
F(b)=F(a)
$$
And thus $F(b)=0$.
\end{proof}

\begin{lemma}\label{K_1=CK_2}Let $\alpha\notin\left\{0, \infty\right\}$. Let $K_1$ be the bilinear form on $V_k(\alpha)\times \mathfrak A$
defined by $K_1(x,y)=F(xy)$. Let $K_2$ be the bilinear form on $V_k(\alpha)\times \mathfrak A$
defined by $K_2(x,y)=F(yx)$. Then there exists an operator $C:V_k(\alpha)\rightarrow V_k(\alpha)$
which has a unique eigenvalue $\alpha$ such that
$$
K_1(x,y)=K_2(Cx,y)
$$
\end{lemma}
\begin{proof}
We proceed by induction on $k$. 

For $k=1$ the operator $C$ is $\alpha \cdot 1$.

Assume that the lemma is true for all $m\leq k$. Consider the case $k+1$.
Let $C_k$ be the operator constructed for the space $V_k(\alpha)$.

Pick a basis in $V_{k+1}(\alpha)$ such that the first $r$ vectors belong to
$V_{k}(\alpha)$. Let $s=\dim V_{k+1}(\alpha)$. For each vector $v_{r+1}$...$v_{s}$
pick an element $a_i\in V_k(\alpha)$ using the definition of the space $V_{k+1}(\alpha)$
with $\alpha_0=0$:
$$
F(v_iy)-\alpha F(yv_i)=F(a_iy)
$$
For $y\in \mathfrak A$ we have
$$
K_1(v_i,y)-\alpha K_2(v_i, y)=K_1(a_i,y)=K_2(C_ka_i, y)
$$
Thus
$$
K_1(v_i,y)=K_2(\alpha v_i+C_ka_i, y)
$$
We now define $C_{k+1}$ by its action on basis vectors of $V_{k+1}(\alpha)$:
$$
C_{k+1}v_i=\left\{
\begin{array}{ll}
1\leq i \leq r \qquad &C_k v_i \\
(r+1)\leq i\leq s \qquad &\alpha v_i + C_k a_i\\
\end{array}
\right.
$$
We see that $C_{k+1}$ has indeed only one eigenvalue $\alpha$.
Also for any basis vector $v_i$ and any $y\in \mathfrak A$ we have
$$
K_1(v_i,y)=K_2(C_{k+1}v_i, y)
$$
Since $\left\{v_i\right\}$ form the basis of $V_{k+1}(\alpha)$ the above equality
holds for any element of $V_{k+1}(\alpha)$.
\end{proof}

\begin{lemma}\label{F(xa)}For $\alpha\neq\infty$ we have
$$
V_{k+1}(\alpha)=\left\{b\in\mathfrak A: \exists a\in V_k(\alpha)\implies\left(\forall x\in\mathfrak A\imply F(bx)-\alpha F(xb)=F(xa)\right)\right\}
$$
In other words, we can set $\alpha_0=\infty$ in the definition of the spaces $V_k(\alpha)$ with
$\alpha\neq \infty$.
\end{lemma}
\begin{proof}
Pick $\alpha_0\neq\alpha$. Let us proceed by induction. For $k=1$ the lemma is true
because the right hand side of the equation in definition of $V_1(\alpha)$ is $0$. 
Assume the lemma holds for all $k\leq n$. From definition of $V_{n+1}(\alpha)$ we have
{\small
$$
\begin{array}{l}
V_{n+1}(\alpha)=\left\{b\in\mathfrak A: \exists a\in V_n(\alpha)\implies\right.\\
\qquad\qquad\qquad\left.\implies\left(\forall x\in\mathfrak A\imply F(bx)-\alpha F(xb)=F(ax)-\alpha_0F(xa)\right)\right\}=\\
\qquad=\left\{b\in\mathfrak A: \exists a\in V_n(\alpha)\exists a'\in V_{n-1}(\alpha)\implies\right.\\
\qquad\qquad\qquad \left.\implies\left(\forall x\in\mathfrak A\imply F(bx)-\alpha F(xb)=\alpha F(xa)+F(xa')-\alpha_0F(xa)\right)\right\}=\\
\qquad=\left\{b\in\mathfrak A: \exists a\in V_n(\alpha)\exists a'\in V_{n-1}(\alpha)\implies\right.\\
\qquad\qquad\qquad \left.\implies\left(\forall x\in\mathfrak A\imply F(bx)-\alpha F(xb)=F\left(x\left((\alpha-\alpha_0) a+a'\right)\right)\right)\right\}
\end{array}
$$
}
We observe that $(\alpha-\alpha_0)a+a'$ is an element of $V_{n}(\alpha)$. 
Thus the lemma holds for $V_{n+1}(\alpha)$ as well.
\end{proof}

{\bf Note.} We observe that theorem \ref{V_k_prop} and lemmas \ref{kernel F}, \ref{K_1=CK_2} and \ref{F(xa)} hold
for any associative algebra $\mathfrak A$.   

\begin{lemma}\label{V_k otimes V_m} Let $\mathfrak A$ and $\mathfrak B$ be two associative algebras. Let $F$ and $G$
be linear functionals on algebras $\mathfrak A$ and $\mathfrak B$ correspondingly. Let $\alpha$ and
$\beta$ be such that $\left\{\alpha,\beta\right\}\neq \left\{0,\infty\right\}$.

Then 
$$
V_k^{\mathfrak A}(\alpha)\otimes V_m^{\mathfrak B}(\beta) \subset V_{k+m-1}^{\mathfrak A\otimes \mathfrak B}(\alpha\beta)
$$
where the latter space was constructed using functional $F\otimes G$.
\end{lemma}
\begin{proof}
First of all, let us note that because the algebra with transposed multiplication law numerates spaces $V_k$
with $1/\alpha$ it is sufficient to prove this lemma in the case of finite $\alpha$ and $\beta$.

Let $b_1\in V_k^{\mathfrak A}(\alpha)$, $b_2\in V_m^{\mathfrak B}(\beta)$, $x\in\mathfrak A$ and $y\in\mathfrak B$.
We have:
{\small
$$
\begin{array}{l}
\left(F\otimes G\right)\left((b_1\otimes b_2)\cdot (x\otimes y)\right):=F(b_1x)G(b_2y)=\\
\qquad\qquad=\left(\alpha F(xb_1)+F(xa_1)\right)\left(\beta G(yb_2)+G(ya_2)\right)=\\
\qquad\qquad=\alpha\beta  F(xb_1)G(yb_2)+\alpha F(xb_1)G(ya_2)+\beta F(xa_1)G(yb_2)+F(xa_1)G(ya_2)=\\
\qquad\qquad=\left(F\otimes G\right)\left((x \otimes y)\cdot (b_1\otimes b_2)+(x \otimes y)\cdot (b_1\otimes a_2+a_1\otimes b_2+ a_1\otimes a_2)\right) \\
\end{array}
$$
}
Therefore the lemma holds for $k=m=1$ as in this case $a_1=a_2=0$. Also the computation above serves
as an induction step in $n=k+m$.
\end{proof}

\begin{lemma}
Let $\mathfrak A$ and $\mathfrak B$ be two associative algebras. Let $F$ and $G$ be two linear functionals
on algebras $\mathfrak A$ and $\mathfrak B$ correspondingly. Then
$$
\Stab_F^{\mathfrak A}(0)\otimes \Stab_G^{\mathfrak B}(\infty)+\Stab_F^{\mathfrak A}(\infty)\otimes \Stab_G^{\mathfrak B}(0)\subset \Nil_{F\otimes G}^{\mathfrak A\otimes \mathfrak B}
$$
\end{lemma}
\begin{proof}
Again because of the argument that algebra with transposed multiplication law numerates spaces $V_k$ with $1/\alpha$
it is sufficient to establish that $\Stab_F^{\mathfrak A}(0)\otimes \Stab_G^{\mathfrak B}(\infty)$ is a subset
of $\Nil_{F\otimes G}^{\mathfrak A\otimes \mathfrak B}$.

Let $b_1\in\Stab_F^{\mathfrak A}(0)$ and $b_2\in\Stab_G^{\mathfrak B}(\infty)$. We have
$$
\left(F\otimes G\right)\left((b_1\otimes b_2)\cdot (x\otimes y)\right):=F(b_1x)G(b_2y)=0\cdot G(b_2y)=0
$$
Also
$$
\left(F\otimes G\right)\left((x\otimes y)\cdot (b_1\otimes b_2)\right):=F(xb_1)G(yb_2)=F(xb_1)\cdot 0=0
$$
Thus $\Stab_F^{\mathfrak A}(0)\otimes \Stab_G^{\mathfrak B}(\infty)\subset \Nil_{F\otimes G}^{\mathfrak A\otimes \mathfrak B}$.
\end{proof}

\begin{lemma}
Let $\mathfrak A$ and $\mathfrak B$ be two associative algebras. Let $F$ and $G$ be two linear functionals
on algebras $\mathfrak A$ and $\mathfrak B$ correspondingly. Then
$$
\begin{array}{ll}
\Stab_F^{\mathfrak A}(0)\otimes\mathfrak B+\mathfrak A\otimes \Stab_G^{\mathfrak B}(0)&\subset \Stab_{F\otimes G}^{\mathfrak A \otimes \mathfrak B}(0) \\
\Stab_F^{\mathfrak A}(\infty)\otimes\mathfrak B+\mathfrak A\otimes \Stab_G^{\mathfrak B}(\infty)&\subset \Stab_{F\otimes G}^{\mathfrak A \otimes \mathfrak B}(\infty) \\
\end{array}
$$
\end{lemma}
\begin{proof}
Consider first the case $\Stab_F^{\mathfrak A}(0)\otimes\mathfrak B\subset \Stab_F^{\mathfrak A \otimes \mathfrak B}(0)$.

Using definition of $\Stab_F^{\mathfrak A}(0)$ we derive:
$$
\begin{array}{l}
\Stab_F^{\mathfrak A}(0)\otimes \mathfrak B = \left\{b_1\in\mathfrak A:\forall x\in\mathfrak A\imply F(b_1x)=0\right\}\otimes\mathfrak B=\\
\qquad\qquad = \span\left( \left\{b_1\otimes b_2: b_1\in\mathfrak A, b_2 \in\mathfrak B, \forall x\in\mathfrak A \forall y\in\mathfrak B\imply \right.\right.\\
\qquad\qquad\qquad\left.\left.\imply\left(F\otimes G\right)((b_1\otimes b_2)(x\otimes y))=F(b_1x)G(b_2y)=0\right\}\right)\subset\\
\qquad\qquad\subset \Stab_{F\otimes G}^{\mathfrak A\otimes \mathfrak B}(0)
\end{array}
$$

Next, observe that by argument of symmetry $\mathfrak A\otimes \mathfrak B \leftrightarrow \mathfrak B\otimes \mathfrak A$
and symmetry $\mathfrak A \leftrightarrow \left(\mathfrak A \textrm { with transposed multiplication} \right)$, 
$\alpha\leftrightarrow 1/\alpha$ we must have as well
$$
\begin{array}{rl}
\mathfrak A\otimes \Stab_G^{\mathfrak B}(0)&\subset \Stab_{F\otimes G}^{\mathfrak A \otimes \mathfrak B}(0) \\
\Stab_F^{\mathfrak A}(\infty)\otimes\mathfrak B&\subset \Stab_{F\otimes G}^{\mathfrak A \otimes \mathfrak B}(\infty) \\
\mathfrak A\otimes \Stab_G^{\mathfrak B}(\infty)&\subset \Stab_{F\otimes G}^{\mathfrak A \otimes \mathfrak B}(\infty) \\
\end{array}
$$
which concludes the proof of this lemma.
\end{proof}

\subsection {Tensor products of type 1 algebras}
In the previous section we have seen that the spaces $V_k(\alpha)$ satisfy some remarkable
properties with respect to tensor products of associative algebras. A natural question
is whether this reflects on the characteristic polynomial of a tensor product of associative
algebras.

\subsubsection{Tensor products of matrices}

Since the definition of characteristic polynomial involves determinant of matrices 
we first turn our attention to tensor products of matrices.

\begin{definition} \buzz{Tensor product of matrices}

Let $A$ and $B$ be two matrices
with coefficients in rings $\mathcal R_1$ and $\mathcal R_2$
respectively. Let commutative ring $\mathcal R$ be a subring of both $\mathcal R_1$
and $\mathcal R_2$.
The tensor product $A \otimes_{\mathcal R} B$ is defined as a block matrix with each block $(i,j)$ having dimensions of matrix $B$
and equal to $A_{i,j}\otimes_{\mathcal R} B$, that is the matrix obtained from $B$ by taking tensor products of
a certain element of $A$ with entrees of $B$. Thus $A \otimes B$ has
coefficients in $\mathcal R_1 \otimes_{\mathcal R} \mathcal R_2$.

\end{definition}

\begin{proposition}The tensor product of matrices has the following
properties:
\begin{enumerate}
\item distributive w.r.t. addition
\item $\left(A\otimes_{\mathcal R}B\right)\cdot \left(C\otimes_{\mathcal R}D\right)=
\left(AC\right)\otimes_{\mathcal R}\left(BD\right)$
\item $\left(A\otimes_{\mathcal R}B\right)^{-1}=\left(A^{-1}\right)\otimes_{\mathcal R}\left(B^{-1}\right)$
\end{enumerate}
\end{proposition}

\begin{theorem}
Let $A$ and $B$ be square matrices of dimensions $k$ and $n$ respectively,
with coefficients in commutative rings $\mathcal R_1$ and $\mathcal R_2$. Let
ring $\mathcal R$ have the property that $\mathcal R \subset \mathcal R_1$
and $\mathcal R \subset \mathcal R_2$. Then
$$
\det\left(A\otimes_{\mathcal R}B\right)=\left(\det
A\right)^n\otimes_{\mathcal R}\left(\det
B\right)^k
$$
\end{theorem}
\begin{proof}
\begin{enumerate}
\item If $A$ and $B$ are diagonal the statement is proved by a simple
computation.
\item Let $R_1=R_2=R=\mathbb{C}$. Let $A=C_1D_1C_1^{-1}$ and
$B=C_2D_2C_2^{-1}$ where $D_1$ and $D_2$ are diagonal. Then 
\begin{eqnarray}
\lefteqn{A\otimes_{\mathcal R}B=\left(C_1D_1C_1^{-1}\right)\otimes_{\mathcal
R}\left(C_2D_2C_2^{-1}\right)= }\nonumber\\
& & \qquad
=\left(C_1\otimes_{\mathcal R}C_2\right)
\left(D_1\otimes_{\mathcal R}D_2\right)
\left(C_1^{-1}\otimes_{\mathcal R}C_2^{-1}\right)=
 \nonumber \\
& & \qquad
=\left(C_1\otimes_{\mathcal R}C_2\right)
\left(D_1\otimes_{\mathcal R}D_2\right)
\left(C_1\otimes_{\mathcal R}C_2\right)^{-1}
 \nonumber 
\end{eqnarray}
and
\begin{eqnarray}
\lefteqn{\det\left(A\otimes_{\mathcal R}B\right)=} \nonumber\\
& &\qquad
=\det\left(\left(C_1\otimes_{\mathcal R}C_2\right)
\left(D_1\otimes_{\mathcal R}D_2\right)
\left(C_1\otimes_{\mathcal R}C_2\right)^{-1} \right)=
\nonumber \\
& &\qquad
=\det\left(
D_1\otimes_{\mathcal R}D_2
 \right)=
\nonumber \\
& & \qquad
=\left(\det D_1\right)^n \left(\det D_2\right)^k=\nonumber\\
& & \qquad
=\left(\det A\right)^n \left(\det B\right)^k\nonumber
\end{eqnarray}
\item Since both sides of the equation $\det\left(A\otimes_{\mathcal R}B\right)=\left(\det
A\right)^n\otimes_{\mathcal R}\left(\det
B\right)^k$ are polynomials in elements of $A$ and $B$ with integral coefficients and we know that over $\mathbb{C}$
all generic $A$ and $B$ satisfy the equation we must have that the polynomials
are identical. This concludes the proof of the theorem.
\end{enumerate}
\end{proof}

\begin{theorem} Let $A$ and $B$ be two matrices with coefficients in commutative
rings $\mathfrak R_1$ and $\mathfrak R_2$ respectively. Let $\mathfrak R$ 
be a unital subring of both $\mathfrak R_1$ and $\mathfrak R_2$. Then there exists
invertible matrix $U$ with coefficients in $\mathfrak R$ that depends only on
dimensions of $A$ and $B$ such that
$$
A\otimes_{\mathfrak R} B=U \left(B\otimes_{\mathfrak R} A\right) U^{-1}
$$
\end{theorem}
\begin{proof}
Let $k$ denote the size of $A$ and $m$ denote the size of $B$. Consider $A$ and 
$B$ as operators acting in $k$-dimensional space $V$ with basis $\left\{f_i\right\}_{i=0}^{k-1}$
and $m$-dimensional space $W$ with basis $\left\{g_i\right\}_{i=0}^{i=m-1}$ correspondingly.

Then the tensor product $A\otimes_{\mathfrak R} B$ is uniquely defined as an operator
in $V\otimes_{\mathfrak R} W$ (which we consider to be a vector space over $\mathfrak R_1\otimes_{\mathfrak R}\mathfrak R_2$).
The matrix representation of operator $A\otimes_{\mathfrak R} B$ depends on the 
choice of basis $\left\{e_i\right\}_{i=0}^{km-1}$ in $V\otimes_{\mathfrak R} W$.

If we use the formula
$$
e_i=f_{[i/m]} \otimes g_{i-[i/m]\cdot m}
$$
we obtain the definition of matrix $A \otimes_{\mathfrak R} B$ (Here $[x]$ denotes the integral part of $x$).

If we use the formula
$$
e_i=f_{i-[i/k]\cdot k} \otimes g_{[i/k]}
$$
we obtain the definition of matrix $B \otimes_{\mathfrak R} A$.

Thus $A \otimes_{\mathfrak R} B$ and $B \otimes_{\mathfrak R} A$ are simply two different matrix representations of the same
operator and the matrix $U$ is the transformation matrix from one basis to the other.
\end{proof}

\begin{theorem}\label{ext_cayley}\buzz{Extended Cayley theorem}
Let $A$ and $B$ be two $n\times n$ matrices over an algebraically closed
field $k$ and $C$ and $D$ be two $m\times m$
matrices over the same field $k$. Define $\chi(\lambda,\mu)=\det(\lambda A+\mu B)$. Then
$$
\det(\lambda A\otimes C +\mu B \otimes D)=det(\chi(\lambda C,\mu D))
$$
\end{theorem}

Before proceeding with the proof we must explain in what sense we 
consider $\chi(\lambda C,\mu D)$. Indeed, matrices $C$ and $D$ might
not commute making $\chi(\lambda C,\mu D)$ ambiguous.
In our situation the right definition is as follows:

First, we notice that $\chi(\lambda,\mu)$ is homogeneous, thus it
can be decomposed into a product of linear forms ($k$ is algebraicly closed):
$$
\chi(\lambda,\mu)=\prod_i(\lambda \alpha_i +\mu \beta_j)
$$

Then we define 
$$
\chi(\lambda C,\mu D)=\prod_i(\lambda \alpha_i C +\mu \beta_j D)
$$

There is still some ambiguity about the order in which we multiply 
linear combinations of $C$ and $D$ but it does not affect the value
of $\det(\chi(\lambda C,\mu D))$.

\begin{proof}

{\bf Step 1.} Let $A$ and $C$ be identity matrices of sizes $n\times n$ and
$m \times m$ respectively. 

Then 
$$
\chi(\lambda,\mu)=\prod_i(\lambda+\mu \gamma_i)
$$
where $\gamma_i$ are eigenvalues of $B$.
\begin{eqnarray}
\det(\chi(\lambda,\mu D))=\det\left(\prod_i(\lambda+\mu \gamma_i D) \right)\qquad\cr
\qquad=\prod_i\det\left(\lambda+\mu \gamma_i D \right)=
\prod_{i,j}\left(\lambda+\mu \gamma_i \epsilon_j\right)
\end{eqnarray}
where $\epsilon_j$ are eigenvalues of $D$.

On the other hand one easily derives that eigenvalues of $B\otimes D$ are 
$\gamma_i \epsilon_j$ and thus
$$
\det(\lambda +\mu B\otimes D)=\prod_{i,j}\left(\lambda+\mu \gamma_i\epsilon_j\right)
=\det(\chi(\lambda,\mu D))
$$

{\bf Step 2.}
Let us assume now only that matrices $A$ and $C$ are invertible.

We have 
$$
\chi(\lambda,\mu)=\det(\lambda A+\mu B)=\det(A)\det(\lambda+\mu A^{-1}B)
$$

Denote $\chi'(\lambda,\mu)=\det(\lambda+\mu A^{-1}B)=\prod_i(\lambda+\mu \gamma_i)$, 
where $\gamma_i$ are eigenvalues of $A^{-1}B$.
\begin{eqnarray}
\lefteqn{\det(\chi(\lambda C,\mu D))=\det\left(\det(A)\prod_i(\lambda C+\mu \gamma_i D)\right)=} \cr
&&=\det(A)^m\prod_i\det(\lambda C+\mu \gamma_i D)=\cr
&&=\det(A)^m\prod_i\left( \det(C)\det(\lambda+\mu \gamma_i C^{-1}D)\right)=\cr
&&=\det(A)^m \det(C)^n \det(\chi'(\lambda,\mu C^{-1}D))\nonumber
\end{eqnarray}

By {\em step 1} we have 
$$
\det(\lambda+\mu (A^{-1}B)\otimes(C^{-1}D))=\det(\chi'(\lambda,\mu C^{-1}D))
$$

Observing also that $\det(A\otimes C)=\det(A)^m\det(C)^n$ we obtain
\begin{eqnarray}
\det(\chi(\lambda C,\mu D))=\det(A\otimes C)\det(\lambda+\mu (A^{-1}B)\otimes(C^{-1}D))=\cr
\qquad =\det(\lambda A\otimes C+\mu B\otimes D)\nonumber
\end{eqnarray}
which is the desired formula.

{\bf Step 3.} Assume now that $\chi(\lambda, \mu)\neq 0$. 
Observe that the right side of the formula
involves only polynomials in entries of matrices $C$ and $D$. Furthermore,
it only involves polynomials in coefficients of $\chi(\lambda, \mu)$ computed
for $C$ and $D$.

We observe that the left
side is polynomial in entries of matrices $A$, $B$, $C$ and $D$. 

A natural question is 
what happens under the symmetry $A\leftrightarrow C$, $B\leftrightarrow D$.

The left part is unchanged:
$$
\begin{array}{l}
\det\left(\lambda A\otimes C + \mu B\otimes D\right)=\det\left(\lambda U\left( C\otimes A \right)U^{-1} + \mu U \left(D\otimes B\right) U^{-1}\right)=\\
\qquad\qquad\qquad=\det\left(\lambda C\otimes A  + \mu D\otimes B \right)
\end{array}
$$

Thus the right hand part is polynomial in coefficients of polynomials $\chi(\lambda, \mu)$
computed for pair $A$ and $B$ and pair $C$ and $D$. This polynomial can be viewed
as another generalization of the resultant of a pair of polynomials.

Since
the restriction that $A$ and $C$ be invertible selects a Zariski open subset,
the formula should hold for all $A$,$B$,$C$,$D$ by continuity.
\end{proof}
\subsubsection{Tensor products of algebras}
We are now in position to prove several results about tensor products of type 1 algebras.

\begin{theorem}\label{chi_of_product} Let $\mathfrak A$ and $\mathfrak B$ be two associative algebras. 
Let $F$ and $G$ be linear functionals on $\mathfrak A$ and $\mathfrak B$ correspondingly.
Let $V_1\subset \mathfrak A$ and $V_2\subset \mathfrak B$ be two linear subspaces.
Then
$$
\chi_{F\otimes G, V_1\otimes V_2}^{\mathfrak A\otimes \mathfrak B}(\lambda, \mu)=\det \chi_{F,V_1}^{\mathfrak A}(\lambda \AFVAVA{B}{G}, \mu \AFVAVA{B}{G}^T)
$$
Here $\AFVAVA{B}{G}$ is the multiplication table of $V_2$ evaluated in $G$.
\end{theorem}
\begin{proof}
Let us pick a basis $\left\{f_i\right\}_{i=0}^{k-1}$ in $V_1$ and a basis
$\left\{g_i\right\}_{i=0}^{m-1}$ in $V_2$. We use $\AF$ to denote the matrix $\left\|F(f_if_j)\right\|$
and $\AFVAVA{B}{G}$ to denote the matrix $\left\|G(g_ig_j)\right\|$.

We define the basis $\left\{e_i\right\}_{i=0}^{km-1}$ in $V_1\otimes V_2$ as following:
$$
e_i=f_{[i/m]}\otimes g_{i-[i/m]\cdot m}
$$
Let $C$ be the matrix $\left\|(F\otimes G)(e_ie_j)\right\|=\left\|F\left(f_{[i/m]}f_{[j/m]}\right)G\left(g_{i-[i/m]\cdot m}g_{j-[j/m]\cdot m}\right)\right\|$. Then
$$
C=\AF\otimes \AFVAVA{B}{G}
$$
Using theorem \ref{ext_cayley} we compute
$$
\begin{array}{l}
\chi_{F\otimes G, V_1\otimes V_2}^{\mathfrak A\otimes \mathfrak B}(\lambda, \mu):=\det\left(\lambda C+\mu C^T\right)=\\
\qquad\qquad\qquad =\det \chi_{F,V_1}^{\mathfrak A}(\lambda \AFVAVA{B}{G}, \mu \AFVAVA{B}{G}^T)
\end{array}
$$
\end{proof}

\begin{theorem}\label{ind_prod} Let $\mathfrak A$ and $\mathfrak B$ be two type 1 associative algebras which are $1$-precise.
Suppose that characteristic polynomial of $\mathfrak A$ has all its roots except $1$ depend non-trivially on $F$.

Then $\ind \mathfrak A\otimes \mathfrak B=\ind \mathfrak A \cdot \ind \mathfrak B$.
\end{theorem}
\begin{proof}
Because characteristic
polynomial of $\mathfrak A$ has all its roots $\alpha\neq 1$ depend non-trivially on $F$ it is possible
to find an open subset of functionals $F$ in which this polynomial is not divisible
by $\lambda$ or $\mu$, i.e. $\Stab_F(0)=\Stab_F(\infty)=\left\{0\right\}$.

For both algebras $\mathfrak A$ and $\mathfrak B$ there exists an open set of 
functionals $F$ that are $1$-regular and $1$-precise. Therefore the spaces $\Stab_F(1)$
are commutative algebras and the symmetric form $Q_F$  is non-degenerate.

Consider the intersection of open sets obtained in first and second paragraphs 
of this proof. Let $F$ be an element of it.

Let $G$ be a regular $1$-precise functional on $\mathfrak B$.

Consider the functional $F\otimes G$ on the algebra $\mathfrak A\otimes \mathfrak B$.
We can compute the value of characteristic polynomial of $\mathfrak A\otimes \mathfrak B$
in $F\otimes G$ by using theorem \ref{chi_of_product}. 

Since the characteristic 
polynomial of $\mathfrak A$ in $F$ is not divisible by either $\lambda$ or $\mu$
we must have that the characteristic polynomial of $\mathfrak A\otimes \mathfrak B$
in $F\otimes G$ is non-zero - and thus the pair $(\mathfrak A\otimes \mathfrak B, F\otimes G)$
is type 1.

All the roots except unity of characteristic polynomial of $\mathfrak A$ depend
non-trivially on $F$. Thus we can vary $F$ to make sure that the only time the 
product $\alpha\beta$ (where $\alpha$ is the root of characteristic polynomial of $\mathfrak A$
in $F$ and $\beta$ is the root of characteristic polynomial of $\mathfrak B$ in $G$) 
is equal to $1$ is when $\alpha=\beta=1$. Therefore, the highest power of $\lambda-\mu$
that divides the characteristic polynomial of $\mathfrak A\otimes \mathfrak B$ in $F\otimes G$
is equal to $\ind \mathfrak A \otimes \ind\mathfrak B$.

By lemma \ref{V_k otimes V_m}, $\Stab_F(1)\otimes \Stab_G(1)\subset \Stab_{F\otimes G}(1)$. Considering
the dimensions it must be $\Stab_F(1)\otimes \Stab_G(1) = \Stab_{F\otimes G}(1)$. Thus $F\otimes G$ is $1$-precise.

Since both $\Stab_F(1)$ and $\Stab_G(1)$ are commutative $\Stab_{F\otimes G}(1)$ is commutative as well. Since symmetric form
$Q_F$ ($Q_G$) on $\Stab_F(1)$ (respectively $\Stab_G(1)$) is non-degenerate it must be that the symmetric form
$Q_{F\otimes G}$ on $\Stab_{F\otimes G}(1)$ is non-degenerate as well. Therefore, by theorem \ref{F is 1-regular}, $F\otimes G$ is $1$-regular.

We have proven that $F\otimes G$ is both $1$-precise and $1$-regular and that $\Stab_{F\otimes G}(1)=\Stab_F(1)\otimes \Stab_G(1)$.
Therefore, $\ind\mathfrak A\otimes \mathfrak B=\ind \mathfrak A \cdot \ind\mathfrak B$.
\end{proof}

{\bf Corollary 1.} The algebra $\mathfrak A\otimes \mathfrak B$ is type 1 and $1$-precise.

{\bf Corollary 2.} Since the algebra of $n\times n$ matrices $\Mat_n$ satisfies conditions on algebra $\mathfrak A$ and
any algebra $\mathfrak B$ with unity and index $1$ satisfies conditions on algebra $\mathfrak B$ we must have
for such algebras
$$
\ind \Mat_n \otimes \mathfrak B=n
$$

\section{Interaction with radical methods}
A classical method of studying associative algebras is to introduce a notion of radical
(see \cite{H,P}). It is reasonable to ask whether the method presented in this paper
provides anything beyound and, in particular, whether Jacobson's radical structure
can be analyzed with this method. Also, one may wonder whether the spectrum $\left\{\alpha(F)\right\}$ 
of the algebra will reflect anything more than the structure of the semisimple factor.

We would like to note that according to the definition in \cite{H} the $P$-radical is
an ideal that is intrinsic to the algebra. However, we know that, for a unital
algebra, the functional $F$ is identically zero on all spaces $\Stab_F(\alpha)$
and $\Nil_F$ except for $\Stab_F(1)$. Furthermore, we know, that for a $1$-regular $F$,
$\Stab_F(1)$ is commutative.

Therefore, if one chooses $1$-regular $F$ so that it does not vanish on at least one element of whatever
radical we are interested in, the radical will either intersect the spaces $\Stab_F(\alpha)$
in a non-trivial way, or be contained entirely within $\Stab_F(1)$ and thus be 
commutative.

We would like to augment this point with the following theorem:
\begin{theorem}
Let $S=\left\{e_i\right\}$ be a countable collection of non-zero elements of $\mathfrak A$.
Then the set of functionals $F$ that does not vanish on all elements of $S$ is 
Baire category $2$ (for definition of Baire categories see \cite{K}).
We consider $\mathfrak A^*$ in Euclidian topology.
\end{theorem}
\begin{proof}
Indeed the set of functionals $F$ that does not vanish on a single element $e_i$
is open and dense in $\mathfrak A^*$. Therefore the set of all functionals $F$
that do not vanish on $S$ is an intersection of open dense sets in a complete
space with Euclidian metric and thus is Baire category $2$.
\end{proof}

Thus, even if we have selected a countable family of non-trivial ideals, we can choose $1$-regular $F$ so that it
does not vanish on any of them.

In regards to the question of whether the spectrum $\left\{\alpha(F)\right\}$ provides anything beyound
characteristics of the semisimple factor, we would like to note that the spectrum
of $\Mat_n$ consists of functions $\alpha(F)$ that, for a generic $F$, depend 
non-trivially on it and do not vanish in any points. As there are plenty of examples
where the spectrum includes constants $\alpha$, in particular $0$ and $\infty$,
this cannot be only due to the semisimple factor.

On the other hand, if we want to discount contribution of a particular ideal $\mathfrak I$
within an algebra, we can either study its factor, or, equivalently, study functionals $F$
pulled back from the factor - this will result in $\mathfrak I\subset \Nil_F$ for all such
functionals.

\section{Open questions}
The identity $\ind \Mat_n \otimes \mathfrak B=n$ does not generalize to any pair of type $1$ algebras as shown
by the following example:

\begin{example} The index of the algebra $UT(2)$ of upper triangular $2\times 2$ matrices is $1$.
The index of $UT(2)\otimes UT(2)$ is $3$.
\end{example}
\begin{proof}
Indeed, the multiplication table of $UT(2)$ is
$$
\begin{array}{ccccc}
&\vline & a & b & c \\
 \hline
a& \vline& a & b & 0 \\
b& \vline& 0 & 0 & b \\
c&\vline& 0 & 0  & c \\
\end{array}
$$
The characteristic polynomial if $UT(2)$ is 
$$
\chi(\lambda, \mu, F)=-\lambda\mu b^2(\lambda+\mu)(a+c)
$$
And thus $\dim\Stab_F(0)=\dim\Stab_F(1)=\dim \Stab_F(\infty)=1$

For the tensor product $UT(2)\otimes UT(2)$ the multiplication table is
$$
\begin{array}{ccccccccccc}
&\vline & a & b & c & d & e & f & g & h & p\\
\hline
a & \vline & a & b & c & d & 0 & 0 & 0 & 0 & 0 \\
b & \vline & 0 & 0 & 0 & 0 & b & d & 0 & 0 & 0 \\
c & \vline & 0 & 0 & 0 & 0 & 0 & 0 & c & d & 0 \\
d & \vline & 0 & 0 & 0 & 0 & 0 & 0 & 0 & 0 & d \\
e & \vline & 0 & 0 & 0 & 0 & e & f & 0 & 0 & 0 \\
f & \vline & 0 & 0 & 0 & 0 & 0 & 0 & 0 & 0 & f \\
g & \vline & 0 & 0 & 0 & 0 & 0 & 0 & g & h & 0 \\
h & \vline & 0 & 0 & 0 & 0 & 0 & 0 & 0 & 0 & h \\
p & \vline & 0 & 0 & 0 & 0 & 0 & 0 & 0 & 0 & p \\
\end{array}
$$

The characteristic polynomial of $UT(2)\otimes UT(2)$ is
$$
\chi(\lambda, \mu, F)=-\lambda^3\mu^3d^3(\lambda+\mu)^3(ch-dg)(fb-ed)(fb+ch+da+dp)
$$

Therefore the only non-zero subspaces $\Stab_F(\alpha)$ are $\Stab_F(0)$, $\Stab_F(1)$ and $\Stab_F(\infty)$
with dimensions that could be anywhere between 1 and 3.

Direct computation yields:
$$
\begin{array}{rl}
\Stab_F(0)=\left\{\right.&(-\frac{d}{a},0,0,1,0,0,0,0,0),\\
	& (-\frac{c}{a}, 0, 1, 0,0,0,0,0,0),\\
	& (-\frac{b}{a}, 1, 0,0,0,0,0,0,0)\left.\right\}\\
\Stab_F(1)=\left\{\right.&(\frac{d}{h}, \frac{f}{h}, 1, -\frac{dch+fbd}{hd^2}, 0, \frac{b}{h},0,\frac{c}{h}, \frac{d}{h}),\\
	&(0, -\frac{f}{d}, 0, \frac{fb}{d^2}, 1, -\frac{b}{d}, 0,0,0), \\
	&(1,\frac{f}{d}, 0, -\frac{fb}{d^2}, 0, \frac{b}{d}, 1, 0,1)\left.\right\}\\
\Stab_F(\infty)=\left\{\right.&(0,0,0,-\frac{p}{d},0,0,0,0,1),\\
   & (0,0,0, -\frac{h}{d}, 0,0,0, 1,0),\\
   & (0,0,0, -\frac{f}{d}, 0,1,0,0,0)\left.\right\}\\
\end{array}
$$
and thus $\ind UT(2)\otimes UT(2)=3$.
\end{proof}

A plausible generalization of the index formula to this case is

{\bf Conjecture:} Let $\mathfrak A$ and $\mathfrak B$ be two type $1$ associative algebras
which resonant spectral values are precise in regular functionals.
Then their product is also a type $1$ algebra and the index is given by this formula:
$$
\ind \mathfrak A\otimes \mathfrak B=\ind \mathfrak A\otimes \mathfrak B+\sum_{1\neq\alpha\in\Spec\mathfrak A\intersect \Spec\mathfrak B} \dim\Stab_{F}^{\mathfrak A}(\alpha)\otimes \dim\Stab_G^{\mathfrak B}(1/\alpha)
$$
Here $\Spec\mathfrak A\intersect \Spec\mathfrak B$ denotes resonant spectral values - namely those constants $\alpha$ for which
$\Stab_F(\alpha)$ and $\Stab_G(\alpha)$ are non-zero for $\alpha$-regular functionals.





\begin{thebibliography}{00}




\bibitem{DK}
{Vladimir Dergachev and Alexandre Kirillov},
{\em Index of Lie algebras of Seaweed type},
{{\em Journal of Lie theory} {\bf 10} No. 2 (2000), 331-343}
\bibitem{D}
{Vladimir Dergachev},
{\em Some properties of index of Lie algebras},
{math.RT/0001042}
\bibitem{DD}
{---},
{\em On associative algebras with unity and Lie index $1$},
{math.RA/0005161}
\bibitem{DDD}
{---},
{\em On a new approach to classification of associative algebras},
{math.RA/0203013}
\bibitem{D1}
{J.Dixmier},
{``Enveloping algebras''},
{American Mathematical Society, Providence (1996) 1-379}
\bibitem{E1}{Elashvili, A. G},
{\em Frobenius Lie algebras},
{Funktsional. Anal. i Prilozhen. {\bf 16} (1982), 94--95}
\bibitem{E2}{---},
{\em On the index of a parabolic subalgebra 
of a semisimple Lie algebra},
{Preprint, 1990}
\bibitem{GK}
{I.M.Gelfand and A.A.Kirillov},
{\em Sur les corps li\'es aux alg\'ebres enveloppantes des alg\'ebres de Lie},
{ Publications math\'ematiques {\bf 31} (1966) 5-20}
\bibitem{H}
{Thomas W.Hugerford},
{Algebra},
{Springer-Verlag, New York, 1974}
\bibitem{K}
{John L.Kelley},
{General topology},
{Springer-Verlag, 1975}
\bibitem{P}
{R.S. Pierce},
{Associative algebras},
{Springer-Verlag, New York, 1982}
\end{thebibliography}
\end{document}